\documentclass[a4paper,11pt]{article}

\usepackage{amssymb}
\usepackage{amsmath}
\usepackage{amsthm}
\usepackage{authblk} % for multiple authors with affiliations
\usepackage{hyperref} % for clickable emails

\usepackage{graphicx}
\usepackage{subcaption}
\usepackage{caption}
\usepackage{color}
\usepackage{tikz}
\usepackage{tikz-3dplot}
\usepackage{bm}
\usepackage{placeins}
\usepackage[margin=2.5cm]{geometry}

\newtheorem{theorem}{Theorem}
\newtheorem{remark}{Remark}
\newtheorem{lemma}{Lemma}

\definecolor{darkgreen}{rgb}{0,0.5,0}
\definecolor{purple}{rgb}{0.75,0,1}
\definecolor{reda}{rgb}{0.8,0,0}
\definecolor{redb}{rgb}{0.6,0,0}
\definecolor{redc}{rgb}{0.4,0,0}
\newcommand{\red}[1]{#1}
\newcommand{\blue}[1]{\textcolor{blue}{#1}}
\newcommand{\purple}[1]{\textcolor{purple}{#1}}

\newcommand{\reda}[1]{\textcolor{reda}{#1}}
\newcommand{\redb}[1]{\textcolor{redb}{#1}}
\newcommand{\redc}[1]{\textcolor{redc}{#1}}

\newcommand{\boldmu}{\bm{\mu}}
\newcommand{\bolda}{\mathbf{a}}
\newcommand{\boldb}{\mathbf{b}}
\newcommand{\boldg}{\mathbf{g}}
\newcommand{\boldn}{\mathbf{n}}
\newcommand{\boldq}{\mathbf{q}}
\newcommand{\boldx}{\mathbf{x}}
\newcommand{\boldv}{\mathbf{v}}
\newcommand{\0}{\mathbf{0}}

\newcommand{\intgammat}{\int_{\Gamma(t)}}
\newcommand{\intgammaht}{\int_{\Gamma_h(t)}}
\newcommand{\intpartialgammat}{\int_{\partial\Gamma(t)}}
\newcommand{\nablagammat}{\nabla_{\Gamma(t)}}
\newcommand{\nablagammaht}{\nabla_{\Gamma_h(t)}}
\newcommand{\deltagammat}{\Delta_{\Gamma(t)}}
\newcommand{\material}{\partial^\bullet}
\newcommand{\materialh}{\partial^\bullet_h}

\title{The evolving surface morphochemical reaction-diffusion system for battery modeling}

\author[1]{Benedetto Bozzini\thanks{\href{mailto:benedetto.bozzini@polimi.it}{benedetto.bozzini@polimi.it}}}
\author[2]{Massimo Frittelli\thanks{Corresponding author: \href{mailto:massimo.frittelli@unisalento.it}{massimo.frittelli@unisalento.it}}}
\author[3,4,5,6]{Anotida Madzvamuse\thanks{\href{mailto:am823@math.ubc.ca}{am823@math.ubc.ca}}}
\author[2]{Ivonne Sgura\thanks{\href{mailto:ivonne.sgura@unisalento.it}{ivonne.sgura@unisalento.it}}}

\affil[1]{Department of Energy, Politecnico di Milano, Via Lambruschini 4, 20156 Milano, Italy}
\affil[2]{Department of Mathematics and Physics ``E. De Giorgi'', University of Salento, Lecce, Italy}
\affil[3]{Mathematics Department, University of British Columbia, 1984 Mathematics Road, Vancouver, V6T 1Z2, British Columbia, Canada}
\affil[4]{Department of Mathematics and Computational Sciences, University of Zimbabwe, Mt Pleasant, Harare, Zimbabwe}
\affil[5]{Department of Mathematics and Applied Mathematics, University of Pretoria, Pretoria, 0132, South Africa}
\affil[6]{Applied Mathematics, University of Johannesburg, PO Box 524, Auckland Park, Johannesburg, 2006, South Africa}

\date{} % remove date

\begin{document}

\maketitle

%% Abstract
\begin{abstract}
%% Text of abstract
It is well known that phase formation by electrodeposition yields films of poorly controllable morphology. This typically leads to a range of technological issues in many fields of electrochemical technology. Presently, a particularly relevant case is that of high-energy density next-generation batteries with metal anodes, that cannot yet reach practical cyclability targets, owing to uncontrolled elelctrode shape evolution. In this scenario, mathematical modelling is a key tool to lay the knowledge-base for materials-science advancements liable to lead to concretely stable battery material architectures. 
In this work, we introduce the Evolving Surface DIB (ESDIB) model, a reaction–diffusion system posed on a dynamically evolving electrode surface. Unlike previous fixed-surface formulations, the ESDIB model couples surface evolution to the local concentration of electrochemical species, allowing the geometry of the electrode itself to adapt in response to deposition. 
To handle the challenges related to the coupling between surface motion and species transport,  we numerically solve the system by proposing an extension of the Lumped Evolving Surface Finite Element Method (LESFEM) for spatial discretisation, combined with an IMEX Euler scheme for time integration.
The model is validated through six numerical experiments, each compared with laboratory images of electrodeposition. Results demonstrate that the ESDIB framework accurately captures branching and dendritic growth, providing a predictive and physically consistent tool for studying metal deposition phenomena in energy storage devices.
\end{abstract}

\section*{Keywords}
Battery modeling; Evolving Surface Reaction-Diffusion Systems (ESRDSs); Turing patterns; Lumped Evolving surface finite element method (LESFEM)

\section*{MSC 2020}
35K57; 35R01; 65M60

\section{Introduction}
The growing global demand for efficient and reliable energy storage systems has spurred intensive research into the physical and chemical mechanisms that govern battery operation. Among various electrochemical technologies, metal-based batteries---such as those employing lithium, sodium, or zinc---offer high theoretical energy densities and have emerged as strong candidates for next-generation storage solutions. However, their performance, safety, and lifetime remain critically limited by morphological instabilities that arise during electrodeposition processes. In particular, the irregular growth of metallic deposits can lead to the formation of dendritic or mossy structures, which increase internal resistance, reduce the active electrode area, and may eventually cause short-circuit failure. These phenomena significantly compromise both the predictability and durability of batteries, motivating the need for accurate mathematical models capable of describing and predicting such complex interfacial growth dynamics. 
At the macroscopic level, such instabilities manifest as capacity fade and premature failure, while at the microscopic level they originate from nonlinear interactions between ionic transport, surface kinetics, and evolving electrode geometry. Capturing this multiscale coupling is therefore crucial for linking the physics of electrodeposition to the long-term reliability of energy storage devices.  Notwithstanding ever-increasing research efforts, that it would be impossible to mention in this context, the topic remains a widely open one and the concertation between experimental studies and mathematical modeling is still largely missing.

In previous works, we introduced a two--species morphochemical reaction--diffusion system referred to as the \textit{DIB} or \textit{morphochemical model},  designed to simulate electrodeposition on flat 2D \cite{bozzini2013spatio} or curved \cite{lacitignola2017turing} cathodic surfaces.  The DIB model, named after its authors, is formulated on a fixed domain or surface and captures the essential coupling between ionic transport and electrochemical reactions during the early stages of deposition.  In the DIB model, one species models surface morphology intended as normal growth velocity \cite{sgura2017xrf},  while the other represents the chemistry.  In a reaction-diffusion system,  spatial structures (patterns) can arise from small perturbations of an equilibrium that is stable in the absence of diffusion, but turns unstable in the presence thereof.  This mechanism is known as Turing instability \cite{turing1952chemical} and the resulting patterns are therefore called Turing patterns.  The DIB model was proven to exhibit, depending on its parameters, several morphological classes of patterns, such as spots, holes,  worms, stripes, labyrinths \cite{bozzini2013spatio}, long worms \cite{frittelli2024mine}, spiral waves \cite{lacitignola2014spatio, lacitignola2019spiral},  and more \cite{dib_cross}. While this framework  inherently assumes a stationary computational domain or surface, the model becomes inadequate once the deposited layer begins to evolve significantly---particularly when branching, self-intersection, or large-scale deformations occur---since these effects fundamentally alter the geometry of the electrode interface. 

In the present work, we extend this framework by introducing a reaction--diffusion system on an evolving surface, which we term the \textit{Evolving Surface DIB (ESDIB)} model. In contrast to the fixed-domain DIB model, the ESDIB model incorporates the electrode surface as a dynamic entity whose shape is part of the \textit{unknown} solution. The surface evolves according to a material velocity that is everywhere normal to the surface itself and depends on the local morphological species governed by the reaction--diffusion equations.  Consequently, surface evolution is intrinsically coupled to the system’s dynamics, allowing the model to naturally describe metal growth, branching, and complex morphological transitions that cannot be captured within a fixed-surface setting.  Following the framework in \cite{barreira2011surface},  we derive the governing equations from a suitable balance law that accounts for mass transport, Fickian diffusion, and Reynolds' transport theorem. Here,  however, this formulation is modified to account for the fact that one of the two species---representing the electrode morphology---is not affected by surface dilation effects.  This evolving-surface framework thus provides a physically consistent and geometrically flexible extension of our earlier approach. 

\red{Several} methods were proposed in the literature for the spatial discretisation of \red{Evolving Surface PDEs} (ESPDEs) with \emph{known} material velocity.  A non-extensive list includes the \red{Evolving Surface Finite Element Method (ESFEM) \cite{dziuk2007finite}},  projected finite elements \cite{tuncer2017projected},  trace finite elements \cite{lehrenfeld2018stabilized, olshanskii2017trace},  kernel methods \cite{wendland2020solving}, and isogeometric analysis \cite{valizadeh2019isogeometric}.  \red{For the spatial discretisation of ESPDEs with \emph{unknown} material velocity,  we mention again the ESFEM \cite{barreira2011surface}.} To numerically approximate the ESDIB system in space, we choose the \textit{Lumped Evolving Surface Finite Element Method (LESFEM)} \cite{frittelli2018numerical},  thanks to its geometric flexibility and implementation simplicity, like in  \cite{frittelli2018numerical}.  We note, however, that in \cite{frittelli2018numerical}  the LESFEM was confined to the case of known isotropic evolution law. Here, instead, the LESFEM is generalised to account for unknown morphology-driven evolution.

We combine the LESFEM spatial discretisation with a modified \textit{implicit-explicit (IMEX) Euler scheme} for time integration. Once again, the modification is tailored to handle the inherent coupling between surface evolution and the unknown species themselves. Therefore,  at each time step,  a new discrete surface and an updated discrete solution need to be computed at once. The combination LESFEM-IMEX Euler preserves the computational efficiency and geometric flexibility shown in our previous works \cite{frittelli2018numerical, frittelli2017lumped, frittelli2019preserving}.  The resulting fully discrete method, \red{which we implemented in MATLAB, } allows robust and accurate simulations of evolving electrode morphologies over extended timescales.  Other time discretisation schemes used in the nonlinear ESPDEs literature include: explicit Runge-Kutta methods with suitable time step restrictions \cite{dziuk2012runge} and BDF methods \cite{kovacs2023maximal}.

Finally, we validate the proposed ESDIB model through six numerical experiments for different choices of the model parameters and corresponding Turing pattern classes, each compared against experimental observations of electrochemical phase formation obtained from laboratory measurements on battery-relevant systems.  The comparisons demonstrate that the evolving-surface formulation captures key qualitative and quantitative features of metallic growth, including the onset of branching and the transition toward dendritic structures, in good agreement with experimental trends.  
It is specially remarkable that the model is capable of capturing a wide range of experimental 3D morphochemical patterns. This impressive property of the DIB model had been highlighted the 2D case in particular in \cite{bozzini2013spatio} for electrodeposition, in \cite{bozzini2011frequency} for dynamic aspects and in \cite{kazemian2022x} for combined precipitation-electrodeposition.
These results highlight the potential of the ESDIB model as a predictive and extensible framework for studying electrodeposition and related interfacial phenomena in electrochemical systems in battery modeling.  In conclusion,  the proposed approach proves that the interplay between mathematical modeling, numerical methods and lab experiments is a key tool to advance the knowledge on material localization processes in batteries in order to enhance their lifespan.  

\subsection*{Outline of the paper}

The remainder of this paper is structured as follows. 
In Section \ref{sec:dib_model}, we introduce the ESDIB model and derive its governing equations from a suitable balance law that accounts for mass transport, Fickian diffusion, and Reynolds' transport theorem. 
In Section \ref{sec:space_discretisation}, we describe the spatial discretisation of the model based on the \textit{Lumped Evolving Surface Finite Element Method (LESFEM)}. 
In Section \ref{sec:time_discretisation}, we present the \textit{IMEX Euler} time discretisation scheme, which is suitably adapted to also approximate surface evolution. 
In Section \ref{sec:numerical_examples}, we report numerical experiments designed to validate the model, comparing simulation results against experimental images obtained from laboratory studies. 
Finally, Section \ref{sec:conclusions} summarises the main conclusions of the work and outlines possible directions for future research.

\section{The DIB model and the derivation of an Evolving Surface \red{DIB} counterpart}
\label{sec:dib_model}

In this section we will recall the original DIB model on stationary flat domains \cite{bozzini2013spatio} and its subsequent \red{Surface DIB (SDIB)} extension to a stationary spherical surface \cite{lacitignola2017turing}. Then, we will further generalise the model by introducing the new Evolving Surface DIB (ESDIB) model.  We will derive the ESDIB model from suitable balance laws, while retaining its morphochemical nature and encompassing the DIB and SDIB models as special cases.

\subsection{The original \red{DIB} model and its interpretation}
The DIB (acronym from the first names of the authors who pioneered it: Deborah Lacitignola, Ivonne Sgura and Benedetto Bozzini) model  of time-space organization in electrochemical phase-formation is an effective and general theoretical platform to model electrodeposition, corrosion and plating-stripping processes in which the morphochemical coupling between electrochemical charge-transfer and evolving surface chemistry controls electrode shape changes. The physico-chemical and mathematical technicalities of the DIB model have been described in a series of papers, among which \cite{bozzini2015weakly, lacitignola2015spatio} are the most representative ones for an introductory reading. To the readers’ benefit, we summarise here the main features of this reaction-diﬀusion  model. Its physical basis is the coupling of morphology and surface chemistry: described with $\eta(\boldx,t)$ and $\theta(\boldx,t)$, respectively, where $\boldx\in \Omega$, with $\Omega\subset\mathbb{R}^2$ is the flat cathodic surface, and $t\in [0,T]$ is the time variable.  Specifically:
\begin{itemize}
\item $\eta(\boldx,t)$, called \emph{morphology}, is the film thickness increment rate;
\item $\theta(\boldx,t) \in [0,1]$, called \emph{chemistry}, is the surface coverage with an electroactive adsorbate.
\end{itemize}
As discussed in \cite[Section 2]{sgura2017xrf}, the physical meaning of $\eta$ entails that the film profile is the graph of the \emph{height  function} $h:\Omega \times [0,T] \rightarrow\mathbb{R}$ defined as follows
\begin{equation}
\label{height_function}
h(\boldx,t) = \int_0^t \eta(\boldx,\tau)\mathrm{d}\tau, \qquad (\boldx,t) \in \Omega\times [0,T],
\end{equation}
see Fig. \ref{fig:dib} for a schematic representation. 
\red{We remark that $h(\boldx,t)$ can exhibit unlimited growth or decay over time,  in agreement with its physical meaning.  In fact,  the cathode can grow (or corrode) without limits as long as there is available space in the battery. }
The DIB model is the following non-dimensionalised system of reaction-diffusion equations  \red{with homogeneous Neumann boundary conditions} that drives the evolution of $\eta$ and $\theta$:
\begin{equation}
\label{dib_model}
\begin{cases}
\dot{\eta} - \Delta \eta = \rho f(\eta,\theta), \qquad (\boldx,t) \in \Omega \times [0,T];\\
\dot{\theta} - d\Delta \theta = \rho g(\eta, \theta), \qquad (\boldx,t) \in \Omega \times [0,T];\\
\red{\nabla \eta \cdot \boldmu = 0, \qquad (\boldx,t) \in \partial\Omega \times [0,T];}\\
\red{\nabla \theta \cdot \boldmu = 0, \qquad (\boldx,t) \in \partial\Omega \times [0,T];}\\
\eta(\cdot,0) = \eta_0,  \qquad \boldx \in \Omega;\\
\theta(\cdot,0) = \theta_0, \qquad \boldx \in \Omega,
\end{cases}
\end{equation}
coupled through the source terms $f$ and $g$, the physical meaning of which we outline here.  \red{In \eqref{dib_model}, $\boldmu: \partial\Omega \rightarrow\mathbb{R}^2$ is the outward unit vector field} and $\rho>0$ rescales the effective domain size,  \cite{frittelli2024matrix}. The source functions for $\eta$ and $\theta$ (kinetics) are given by:
\begin{align}
\label{f}
&f(\eta, \theta) = A_1(1-\theta)\eta  - A_2\eta^3 - B(\theta-\alpha);\\
\label{g}
&g(\eta,\theta) = C(1+k_2\eta)(1-\theta)[1-\gamma(1-\theta)]-D(1+k_3\eta) \theta(1+\gamma\theta).
\end{align}
with all constants real and positive or equal to zero.  In \eqref{f}-\eqref{g}:
\begin{enumerate}
\item $A_1(1-\theta)\eta$ accounts for the charge-transfer rate at electrode sites free from adsorbate;
\item $-A_2\eta^3$ models mass-transport control;
\item $-B(\theta-\alpha)$ control describes the effect of adsorbates on the growth rate;
\item $C$ and $D$ scale the effect of surface chemistry and local electrochemical activity on adsorption and desorption, respectively;
\item $\gamma$,  $k_2$ and $k_3$ gauge the impact of chemical and electrochemical contributions on adsorbate film formation.
\end{enumerate}
\red{As shown in \cite{bozzini2013spatio}, the spatially homogeneous equilibrium
\begin{equation}
\label{equilibrium_dib}
(\eta_{\text{eq}}, \theta_{\text{eq}}) := (0,\alpha),
\end{equation}
is subject to \emph{diffusion-driven} instability if the model parameters are chosen in a suitable subset of the parameter space, known as \emph{Turing space}.  This means that, by choosing the model parameters in the Turing space and the initial condition $(\eta_0(\boldx),\theta_0(\boldx))$ of \eqref{dib_model} as a small random perturbation of the equilibrium \eqref{equilibrium_dib}, the solution $(\eta(\boldx),\theta(\boldx))$ reaches asymptotically in time a steady state $(\overline{\eta}(\boldx), \overline{\theta}(\boldx))$ known as a \emph{Turing pattern}.  The DIB model was proven to generate several classes of Turing patterns depending on the model parameters,  such as spots, stripes, labyrinths, worms, holes, and more.} Extensive numerical and physico-mathematical analyses have cogently corroborated the capability and ﬂexibility of DIB in modelling for the vast corpus of patterns found in the experimental electrochemical literature, including electrodeposition, corrosion, electrocatalysis and battery science and technology, described in over 30 journal papers, spanning the period 2007-2025.

\subsection{The \red{SDIB} model on a stationary surface}
To encompass the more realistic case when the cathodic surface is not necessarily flat,  the \red{Surface DIB (SDIB)} model on stationary surfaces,  introduced in \cite{lacitignola2017turing}, extends the DIB model \eqref{dib_model}.  Given a compact surface $\Gamma \subset\mathbb{R}^3$ without boundary, the SDIB model is the non-dimensionalised reaction-diffusion system:
\begin{equation}
\label{sdib_model}
\begin{cases}
\dot{\eta} - \Delta_\Gamma \eta = \rho f(\eta,\theta), \qquad (\boldx,t) \in \Gamma \times [0,T];\\
\dot{\theta} - d\Delta_\Gamma \theta = \rho g(\eta, \theta), \qquad (\boldx,t) \in \Gamma \times [0,T];\\
\eta(\cdot,0) = \eta_0,  \qquad \boldx \in \Gamma;\\
\theta(\cdot,0) = \theta_0, \qquad \boldx \in \Gamma,
\end{cases}
\end{equation}
where $\Delta_\Gamma$ is the Laplace-Beltrami operator on $\Gamma$.  \red{Similarly with the DIB model \eqref{dib_model}, the equilibrium \eqref{equilibrium_dib} is subject to diffusion-driven instability and retains the same Turing space,  see \cite{lacitignola2017turing}.  The difference, compared to the DIB model \eqref{dib_model},  is that the excitable modes depend on the shape and size of the surface $\Gamma$.} In analogy with the DIB model \eqref{dib_model},  $\eta(\boldx,t)$ (morphology) is the film thickness increment rate \emph{in the normal direction} w.r.t.  $\Gamma$. Therefore, the term \emph{height function} is no longer appropriate.  For the sake of clarity,  in this context where $\Gamma$ is a general surface, the height function $h(x,t)$ defined in \eqref{height_function} is generalised by the \emph{thickness function} $h_\Gamma: \Gamma \times [0,T] \rightarrow\mathbb{R}$ defined as
\begin{equation}
\label{thickness_function}
h_\Gamma(\boldx,t) = \int_0^t \eta(\boldx,\tau)\mathrm{d}\tau, \qquad (\boldx,t) \in \Gamma\times [0,T].
\end{equation}
\red{Similarly to the height function \eqref{height_function} considered in the flat case,  the thickness function $h(\boldx,t)$ can grow or decay in time without limits.}
Since the spatial domain of $h_\Gamma$ is a possibly curved surface $\Gamma$, we need to generalise the notion of graph of a function to $\Gamma$-variate functions.  Specifically, the film profile is the graph of $h_\Gamma$ defined by
\begin{equation}
\label{SDIB_auxiliary_surface}
\Gamma_N(t) = \left\{\boldx_N(t) := \boldx + h(\boldx,t)\boldn(\boldx) \ \middle|\  \boldx\in\Gamma \right\},
\end{equation}
where $\boldn(\cdot)$ is the outward unit normal vector field on the surface $\Gamma$.  A pictorial representation is given in Figure \ref{fig:sdib}.  This implies that the DIB model successfully models metal growth as long as the profile can be represented as the graph of a function posed on the initial surface.  Clearly, this ceases to be true on the onset of branching.  In conclusion, the SDIB model \eqref{sdib_model} relies on the following simplifying assumption:
\begin{equation}
\label{assumption_stationary_domain}
\begin{split}
&\text{The evolution of the film profile $\Gamma_N(t)$ has a negligible effect on the dynamics and is}\\
&\text{therefore not accounted for in the governing equations}.
\end{split}
\end{equation}
The ESDIB model, which we introduce in the next section, is designed to overcome this limitation.

\subsection{Evolving surface reaction-diffusion systems: background}
\begin{figure}[ht!]
\centering
\begin{subfigure}{\textwidth}
\centering
\begin{tikzpicture}[scale = 0.5]
    % Initial curve (circle)
    \draw[thick, red] (-7,0) -- (7,0) node[right] {$\Omega$};
    
    % Intermediate curve
    \draw[thick, purple, domain=-7:7, smooth, samples=50] plot (
        {\x}, 
        {1+0.5*sin(100*\x)}
    ) node[right] {$y = h(x,t_1)$};
    
    % Outer curve
    \draw[thick, blue, domain=-7:7, smooth, samples=50] plot (
        {\x}, 
        {2+sin(100*\x)} 
    ) node[right] {$y = h(x,t_2)$};
    
    % Black points on the curves
    \filldraw[black] (0.95,0) circle(4pt) node[below] {$\mathbf{x}$};
    \filldraw[black] (0.95,1.5) circle(4pt) node[right] {};
    \filldraw[black] (0.95,3) circle(4pt) node[right] {};
    
    % Draw  arrows
    \draw[thick, ->] (0.95,0) -- (0.95,1.4);
    \draw[thick, ->] (0.95,1.5) -- (0.95,2.9);
    
\end{tikzpicture}
\caption{In the DIB model \eqref{dib_model},  posed on the stationary flat domain $\Omega$, the evolution of the electrode shape is modeled by the graph of the height function $h$ defined in \eqref{height_function}. The domain $\Omega$ remains fixed, while the graph of $h$ evolves in the direction normal to $\Omega$.}
\label{fig:dib}
\end{subfigure}
\begin{subfigure}{0.48\textwidth}
\centering
\begin{tikzpicture}[scale = 0.5]
    % Initial curve (circle)
    \draw[thick, red] (0,0) circle(2) node[above] {$\Gamma$};
    
    % Intermediate curve
    \draw[thick, purple, domain=0:360, smooth, samples=50] plot (
        {(3.5 + 0.7*sin(5*\x))*sin(\x)}, 
        {(3.5 + 0.7*sin(5*\x))*cos(\x)}
    );
    
    % Outer curve
    \draw[thick, blue, domain=0:360, smooth, samples=50] plot (
        {(4.7 + 1*sin(5*\x))*sin(\x)}, 
        {(4.7 + 1*sin(5*\x))*cos(\x)} 
    );
    
    % Black points on the curves
    \filldraw[black] (2,0.25) circle(4pt) node[left] {$\mathbf{x}$};
    \filldraw[black] (4.15,0.3) circle(4pt) node[right] {};
    \filldraw[black] (5.65,0.35) circle(4pt) node[right] {};
    
    % Draw  arrows
    \draw[thick, ->] (2.2,0.25) -- (3.67,0.29);
    \draw[thick, ->] (4.15,0.3) -- (5.4,0.34);
    
    % Labels with dotted lines
    \draw[dotted, thick] (4.15,0.4) -- (4.15,-3);
    \node[below] at (4.15,-3) {$\boldx_N(t_1)$};
    
    \draw[dotted, thick] (5.65,0.35) -- (5.65,-2);
    \node[below] at (5.65,-2) {$\boldx_N(t_2)$};
    
    \draw[dotted, thick] (-3,2.8) -- (-3,4);
    \node[above] at (-3,4) {$\purple{\Gamma_N(t_1)}$};
    
    \draw[dotted, thick] (-1.5,2.7) -- (-1.5,5);
    \node[above] at (-1.5,5) {$\blue{\Gamma_N(t_2)}$};
    
    % Normal direction to intermediate curve
    \draw[thick, purple,->] (4.15,0.3) -- (6.5,1.7) node[above] {$\boldn(\boldx_N(t_1))$};
\end{tikzpicture}
\caption{In the SDIB model \eqref{sdib_model},  posed on the stationary surface $\Gamma$, the evolution of the electrode shape is modeled by the graph $\Gamma_N(t)$,  defined in \eqref{SDIB_auxiliary_surface}, of the thickness function $h_\Gamma$ defined in \eqref{thickness_function}. The domain $\Gamma$ remains fixed, while the graph $\Gamma_N(t)$ evolves in the direction normal to $\Gamma$, not to $\Gamma_N(t)$ itself.}
\label{fig:sdib}
\end{subfigure}
\hfill
\begin{subfigure}{0.48\textwidth}
\centering
\begin{tikzpicture}[scale = 0.5]
    % Initial curve (circle)
    \draw[thick, red] (0,0) circle(2) node[above] {$\Gamma(t_0)$};
    
    % Intermediate curve
    \draw[thick, reda, domain=0:360, smooth, samples=100] plot (
        {(3)*(cos(\x) + 0.12*cos(6*\x))}, 
        {(3)*(sin(\x) + 0.12*sin(6*\x))}
    );
    
    % Outer curve
    \draw[thick, redb, domain=0:360, smooth, samples=100] plot (
        {(3.7)*(cos(\x) + 0.16*cos(6*\x))}, 
        {(3.7)*(sin(\x) + 0.16*sin(6*\x))}
    );
    
    % Outer curve
    \draw[thick, redc, domain=0:360, smooth, samples=100] plot (
        {(4.5)*(cos(\x) + 0.23*cos(6*\x))}, 
        {(4.5)*(sin(\x) + 0.23*sin(6*\x))}
    );
    
    % Black points on the curves
    \filldraw[black] (2,0.25) circle(4pt) node[left] {$\mathbf{x}$};
    \filldraw[black] (3.25,0.43) circle(4pt) node[right] {};
    \filldraw[black] (4.05,0.95) circle(4pt) node[right] {};
    \filldraw[black] (4.6,2) circle(4pt) node[right] {};
    
    % Draw  arrows
    \draw[thick, ->] (2.2,0.25) to[out =5, in = 200] (3.1,0.4);
    \draw[thick, ->] (3.3,0.44) to[out = 20, in = 220] (3.95,0.87);
    \draw[thick, ->] (4.15,1.05) to[out = 45, in = 255] (4.55,1.9);
    
    % Labels with dotted lines
    \draw[dotted, thick] (3.25,0.4) -- (4.5,-3);
    \node[below] at (4.5,-3) {$\boldx(t_1)$};
    
    \draw[dotted, thick] (4.05,0.95) -- (5.2,-2);
    \node[below right] at (5.2,-2) {$\boldx(t_2)$};
    
    \draw[dotted, thick] (4.6,2) -- (6,-1);
    \node[below right] at (6,-1) {$\boldx(t_3)$};
    
    \draw[dotted, thick] (-1.5,2.7) -- (-1.5,5);
    \node[above] at (-1.5,5) {$\reda{\Gamma(t_1)}$};
    
    \draw[dotted, thick] (-2.8,3.2) -- (-2.8,4.2);
    \node[above] at (-2.8,4.2) {$\redb{\Gamma(t_2)}$};
    
    \draw[dotted, thick] (-4,4) -- (-4,5);
    \node[above] at (-4,5) {$\redc{\Gamma(t_3)}$};
\end{tikzpicture}
\caption{In the ESDIB model \eqref{model_evolving}, the electrode profile is given by the evolving surface $\Gamma(t)$, which evolves according to the law \eqref{evolution_law}.  The ESDIB model is posed on the evolving surface $\Gamma(t)$, which evolves normally to $\Gamma(t)$. Therefore,  $\Gamma(t)$ is itself an unknown of the problem. \red{The ESDIB model is physically meaningful up to the onset of self-intersections in $\Gamma(t)$, but is mathematically well-posed at all times.}}
\label{fig:esdib}
\end{subfigure}
\caption{How the evolution of the electrode shape is modeled (a) in the DIB model \eqref{dib_model}, (b) in the SDIB model \eqref{sdib_model} and (c) in the ESDIB model \eqref{model_evolving}.}
\label{fig:concentration_driven_evolution}
\end{figure}
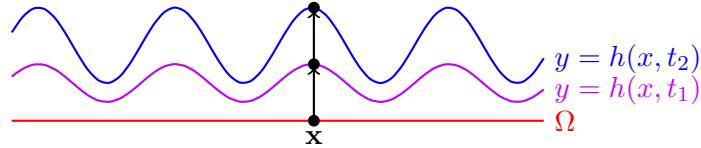
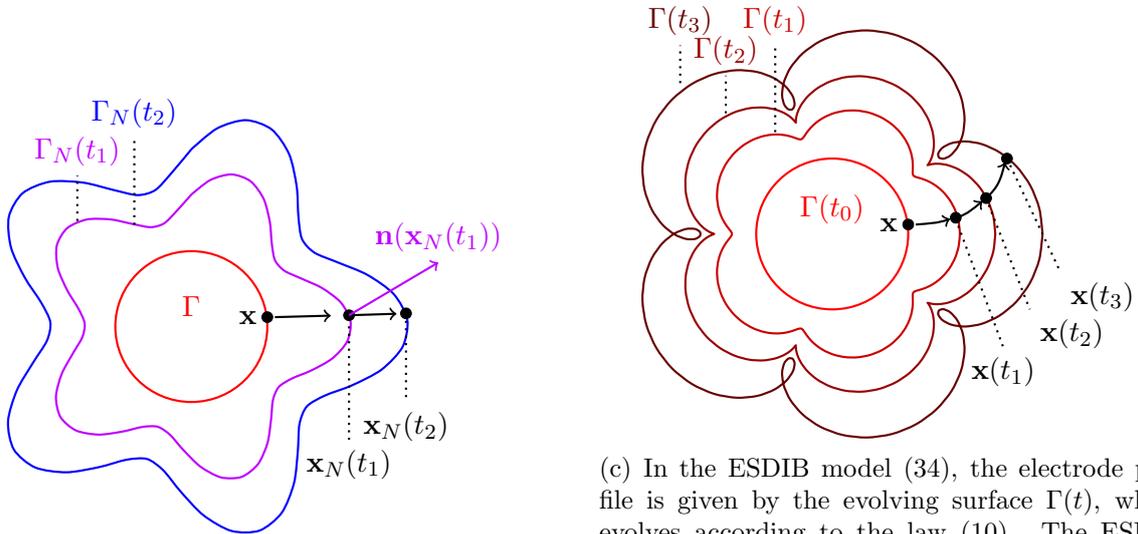
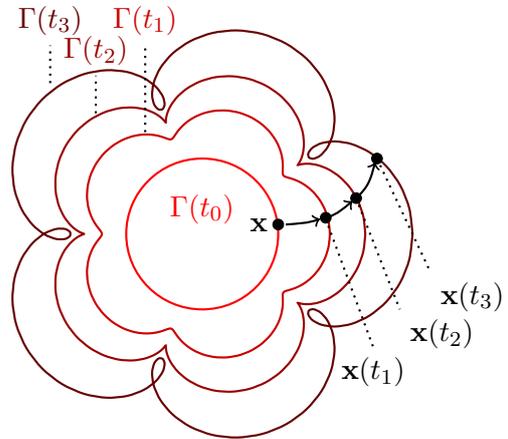

If we remove Assumption \eqref{assumption_stationary_domain},  the spatial domain becomes a time-dependent surface $\Gamma(t)$ \red{originating} from an initial surface $\Gamma_0$.  Thanks to the physical meaning of $\eta$, each point $\boldx$ of $\Gamma(t)$ evolves according to the law
\begin{equation}
\label{evolution_law}
\begin{cases}
\dot{\boldx}(t) := \kappa \eta(\boldx, t) \boldn(\boldx, t), \qquad t\in [0,T];\\
\boldx(0) = \boldx_0 \in \Gamma_0,
\end{cases}
\end{equation}
where $\boldn(\cdot,t)$ is the outward unit normal vector field on the surface $\Gamma(t)$ and $\kappa > 0$ is a suitable proportionality constant.  A pictorial representation is given in Figure \ref{fig:esdib}.  \red{We remark that the evolution law \eqref{evolution_law} does not imply volume conservation and, experimentally,  can lead to exponential growth of surface area. This behavior is physically correct, as discussed in Section \ref{sec:numerical_examples}}. The flow of the Cauchy problem \eqref{evolution_law} is the function $G:\Gamma_0 \times [0,T] \rightarrow\mathbb{R}^3$ such that, for any $\boldx_0\in\Gamma_0$,  $\boldx(t) = G(\boldx_0,t)$ is the solution to the Cauchy problem \eqref{evolution_law}. The evolution law \eqref{evolution_law} is a special case of a general evolution law of the type
\begin{equation}
\label{general_evolution_law}
\dot{\boldx}(t) = \boldv(\eta(\boldx,t),\theta(\boldx,t),\boldx,t),   \qquad t\in [0,T], \boldx\in\Gamma(t).
\end{equation}
Surface PDEs where the spatial domain $\Gamma(t)$ evolves according to a law of the type \eqref{general_evolution_law} are called \emph{evolving surface PDEs (ESPDEs) with concentration-driven evolution}. 
\red{ESPDEs and their derivation from balance laws were systematically introduced in \cite{barreira2011surface}},  see also \cite{dziuk2013finite} for a comprehensive review on the topic.  In the following, to simplify the notations, we recall from \cite{frittelli2018numerical} the definition of \emph{graph of the evolving surface}:
\begin{equation}
\mathcal{G} =: \bigcup_{t\in [0,T]} \Gamma(t) \times \{t\}.
\end{equation}
\red{Similarly, we define the \emph{graph of the boundary of the evolving surface}:
\begin{equation}
\mathcal{G}_\partial := \bigcup_{t\in [0,T]} \partial\Gamma(t) \times \{t\}.
\end{equation}}
The statement $t\in [0,T]$, $\boldx\in \Gamma(t)$ is thus equivalent to $(\boldx,t) \in \mathcal{G}$, \red{and $t\in [0,T]$, $\boldx\in \partial\Gamma(t)$ is equivalent to $(\boldx,t) \in \mathcal{G}_\partial$.} Given any sufficiently smooth function $u:\mathcal{G} \rightarrow\mathbb{R}$, the \emph{material derivative} \cite{frittelli2018numerical} defined by
\begin{equation}
\label{material_derivative}
\partial^\bullet u(\boldx,t) := \frac{\mathrm{d}}{\mathrm{d}t} u(\boldx(t),t), \qquad (\boldx,t) \in \mathcal{G},
\end{equation}
where $\boldx(t)$ solves \eqref{evolution_law}. The material derivative \eqref{material_derivative} can be viewed as a time derivative in a moving reference system. An equivalent way of expressing the material derivative is
\begin{equation}
\label{material_derivative_extension}
\partial^\bullet u(\boldx,t) := \frac{\partial \widetilde{u}}{\partial t} + \boldv\cdot \nabla \widetilde{u}, \qquad t \in [0,T], \ \boldx\in \mathcal{N}(t),
\end{equation}
where $\widetilde{u}:\mathcal{N}(t)\rightarrow\mathbb{R}$ is any smooth extension of $u$ defined on an open neighborhood $\mathcal{N}(t)$ of $\Gamma(t)$, see \cite{dziuk2013finite}.  We now recall some preliminary results.

\begin{theorem}[Integration by parts on evolving surfaces, \cite{dziuk2013finite}]
If $\boldg:\mathcal{G} \rightarrow\mathbb{R}^3$ is a sufficiently smooth vector field tangent to $\Gamma(t)$ at all times, it holds that
\begin{equation}
\label{integration_by_parts}
\intpartialgammat \boldg\cdot\boldmu = \intgammat \nablagammat\cdot \boldg, \qquad t\in [0,T],
\end{equation}
where $\boldmu:\partial \Gamma(t) \rightarrow\mathbb{R}^3$ is the unit outward \emph{conormal} vector field on $\partial \Gamma(t)$,  see \cite{dziuk2013finite}.
\end{theorem}

\begin{theorem}[Reynolds transport theorem, \cite{dziuk2013finite}]
If $g:\mathcal{G} \rightarrow\mathbb{R}$ is a sufficiently smooth scalar function, it holds that
\begin{equation}
\label{reynolds}
\frac{\mathrm{d}}{\mathrm{d}t} \intgammat g = \intgammat (\material g + g\nablagammat\cdot\dot{\boldx}), \qquad t\in [0,T].
\end{equation}
\end{theorem}

\begin{theorem}[Green's formula on surfaces,  \cite{dziuk2013finite}]
If $f,g:\mathcal{G} \rightarrow\mathbb{R}$ are sufficiently smooth scalar functions, it holds that
\begin{equation}
\label{green}
\intgammat \nablagammat f \cdot \nablagammat g = \intgammat f\deltagammat g + \intpartialgammat f\nablagammat g \cdot \boldmu, \qquad t\in [0,T].
\end{equation}
\end{theorem}

\noindent
Let $\mathcal{R}_0 \subset \Gamma_0$ be an arbitrary smooth portion of the initial surface $\Gamma_0$ and let $\mathcal{R}(t) := G(\mathcal{R}_0, t)$ be its corresponding evolved counterpart, where $G$ is the flow defined below \eqref{evolution_law}. As shown in \cite[Section 2.2]{frittelli2018numerical}, the desired governing equations are derived from the balance laws
\begin{align}
\label{balance_law_eta}
&\frac{\mathrm{d}}{\mathrm{d}t} \int_{\mathcal{R}(t)} \eta(\boldx,t)\mathrm{d}\boldx = -\int_{\partial\mathcal{R}(t)} \boldq_\eta \cdot \boldmu +  \int_{\mathcal{R}(t)} \rho f(\eta,\theta)\mathrm{d}\boldx, \qquad t\in [0,T];\\
\label{balance_law_theta}
&\frac{\mathrm{d}}{\mathrm{d}t} \int_{\mathcal{R}(t)} \theta(\boldx,t)\mathrm{d}\boldx = -\int_{\partial\mathcal{R}(t)} \boldq_\theta \cdot \boldmu + \int_{\mathcal{R}(t)} \rho g(\eta,\theta)\mathrm{d}\boldx, \qquad t\in [0,T],
\end{align}
where $\boldmu:\partial \mathcal{R}(t) \rightarrow\mathbb{R}^3$ is the unit outward conormal vector field, this time on $\partial \mathcal{R}(t)$,  see \cite{dziuk2013finite},  $\boldq_\eta$ and $\boldq_\theta$ are the outward fluxes of $\eta$ and $\theta$, respectively, across $\partial \mathcal{R}(t)$, $f$ and $g$ are sufficiently regular reaction kinetics and $\rho>0$ is a rescaling parameter.
\begin{remark}[Interpretation of the balance laws \eqref{balance_law_eta}-\eqref{balance_law_theta}]
To give an empirical interpretation of \eqref{balance_law_eta}-\eqref{balance_law_theta},  we will consider special case $\boldq_\eta = \0$ and  $f=g=0$. In this case,  \eqref{balance_law_eta}-\eqref{balance_law_theta} become
\begin{align}
&\frac{\mathrm{d}}{\mathrm{d}t} \int_{\Gamma(t)} \eta(\boldx,t)\mathrm{d}\boldx =0, \qquad \text{and} \qquad \frac{\mathrm{d}}{\mathrm{d}t} \int_{\Gamma(t)} \theta(\boldx,t)\mathrm{d}\boldx = 0, \qquad t\in [0,T],
\end{align}
i.e.  the total amounts of $\eta$ and $\theta$ are preserved. This means that, in the  case $\partial\Gamma(t) = \emptyset$ and $f=g=0$,  the terms $\eta \nabla_{\Gamma(t)} \cdot \dot{\boldx}(t)$ and $\theta \nabla_{\Gamma(t)} \cdot \dot{\boldx}(t)$ in \eqref{esrds} can be viewed as \emph{dilution terms}: they cause the local concentration of each species to increase or decrease in such a way that their spatial integrals are preserved over time,  see Fig. \ref{fig:type_1_evolution} for an illustration.  
\end{remark}

\noindent
It was proven in \cite{barreira2011surface} that the evolution of two species $\eta(\boldx,t),\theta(\boldx,t)$ fulfilling the balance laws \eqref{balance_law_eta}-\eqref{balance_law_theta} and Fickian diffusion
\begin{equation}
\label{fick}
\boldq_\eta = -\nablagammat \eta, \qquad \text{and} \qquad \boldq_\theta = -d \nablagammat \theta, \qquad (\boldx,t) \in \mathcal{G},
\end{equation}
is governed by the following \red{non-dimensional} \emph{evolving surface reaction-diffusion system (ESRDS)}:
\begin{equation}
\label{esrds}
\begin{cases}
\partial^\bullet \eta + \eta \nabla_{\Gamma(t)} \cdot \dot{\boldx}(t) - \Delta_{\Gamma(t)} \eta = \rho f(\eta,\theta), \qquad (\boldx,t) \in\mathcal{G};\\
\partial^\bullet \theta + \theta \nabla_{\Gamma(t)} \cdot \dot{\boldx}(t) - d\Delta_{\Gamma(t)} \theta = \rho g(\eta, \theta),\qquad (\boldx,t) \in\mathcal{G};\\
\red{\nablagammat\eta \cdot \boldmu = 0, \qquad (\boldx,t) \in\mathcal{G}_\partial;}\\
\red{\nablagammat\theta \cdot \boldmu = 0, \qquad (\boldx,t) \in\mathcal{G}_\partial;}\\
\dot{\boldx} = \boldv(\eta,\theta,\boldx,t), \qquad (\boldx,t) \in\mathcal{G};\\
\eta(\cdot,0) = \eta_0,  \qquad \boldx \in \Gamma(0);\\
\theta(\cdot,0) = \theta_0, \qquad \boldx \in \Gamma(0);\\
\Gamma(0) = \Gamma_0.
\end{cases}
\end{equation}
Due to its clear interpretation, the ESRDS \eqref{esrds} is a general platform for modeling diverse physical and chemical phenomena \cite{barreira2011surface, lacitignola2022pattern, lefevre2010reaction},  and was therefore extensively studied and analysed in the literature \cite{tuncer2017projected, frittelli2018numerical,dziuk2013finite}. The derivation of the ESRDS \eqref{esrds} relies on integration by parts \eqref{integration_by_parts}, Reynold's transport formula \eqref{reynolds}, and Green's formula on surfaces \eqref{green}, see \cite{barreira2011surface}.  In the next section, starting from suitably modified balance laws,  we will derive an evolving surface counterpart of the SDIB model \eqref{sdib_model} using the same tools.

\subsection{Introducing the ESDIB model}
We \red{proceed} to distinguish two kinds of surface evolution.  The surface can (i) be \emph{transported by the material velocity} or (ii) can \emph{evolve by addition or removal of material onto the pre-existing surface}.  The evolution of the same surface $\Gamma(t)$ can be of type (i) or type (ii) for each species, according to the physical meaning thereof.  A species $\theta$ that drives type (i) evolution fulfils a balance law of the type \eqref{balance_law_theta}.  A species $\eta$ that drives type (ii) evolution,  instead,  fulfils the following modified balance law instead of \eqref{balance_law_eta}:
\begin{equation}
\label{balance_law_eta_type2}
\frac{\mathrm{d}}{\mathrm{d}t} \int_{\mathcal{R}(t)} \eta(\boldx,t)\mathrm{d}\boldx = -\int_{\partial\mathcal{R}(t)} \boldq_\eta \cdot \boldmu + \int_{\mathcal{R}(t)} \rho f(\eta,\theta)\mathrm{d}\boldx + \int_{\mathcal{R}(t)} \eta\nabla_{\Gamma(t)} \cdot\dot{\boldx}(t)\mathrm{d}\boldx,  \qquad t\in [0,T].
\end{equation}

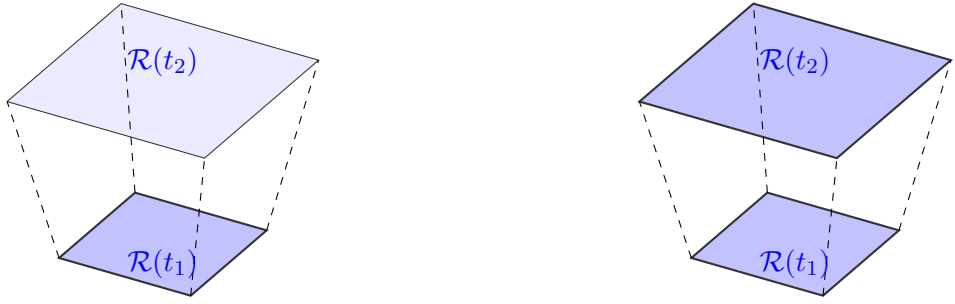
\begin{figure}
\begin{subfigure}{0.48\textwidth}
\tdplotsetmaincoords{60}{120}  % Set view angle
\centering
\begin{tikzpicture}[tdplot_main_coords]

% Define height of the original blue square
\def\zlevel{1.5}
% Define vertical offset for shadow (above the square)
\def\shadowz{4}
% Define scale factor for shadow (bigger than the square)
\def\scale{1.5}

% Original blue square at z = \zlevel
\coordinate (A) at (0,0,\zlevel);
\coordinate (B) at (2,0,\zlevel);
\coordinate (C) at (2,2,\zlevel);
\coordinate (D) at (0,2,\zlevel);

% Shadow projected vertically upward and scaled from center (1,1)
\coordinate (As) at ({1 + \scale*(-1)}, {1 + \scale*(-1)}, \shadowz);
\coordinate (Bs) at ({1 + \scale*(1)},  {1 + \scale*(-1)}, \shadowz);
\coordinate (Cs) at ({1 + \scale*(1)},  {1 + \scale*(1)},  \shadowz);
\coordinate (Ds) at ({1 + \scale*(-1)}, {1 + \scale*(1)},  \shadowz);

% Draw the original blue square
\draw[thick,fill=blue!30,opacity=0.8] (A) -- (B) -- (C) -- (D) -- cycle;

% Draw the lighter blue shadow above
\draw[fill=blue!10,opacity=0.8] (As) -- (Bs) -- (Cs) -- (Ds) -- cycle;

% Draw vertical projection lines
\draw[dashed] (A) -- (As);
\draw[dashed] (B) -- (Bs);
\draw[dashed] (C) -- (Cs);
\draw[dashed] (D) -- (Ds);

% Add labels in blue font
\node[blue] at (1,1,\zlevel - 0.3) {$\mathcal{R}(t_1)$};
\node[blue] at (1,1,\shadowz + 0.3) {$\mathcal{R}(t_2)$};
\end{tikzpicture}
\caption{Type (i) evolution in the absence of reaction and boundary fluxes: dilution ensures that the \emph{spatial integral} on $\mathcal{R}(t)$ of a spatially uniform species is preserved over time.}
\label{fig:type_1_evolution}
\end{subfigure}
\hfill
\begin{subfigure}{0.48\textwidth}
\tdplotsetmaincoords{60}{120}  % Set view angle
\centering
\begin{tikzpicture}[tdplot_main_coords]

% Define height of the original blue square
\def\zlevel{1.5}
% Define vertical offset for shadow (above the square)
\def\shadowz{4}
% Define scale factor for shadow (bigger than the square)
\def\scale{1.5}

% Original blue square at z = \zlevel
\coordinate (A) at (0,0,\zlevel);
\coordinate (B) at (2,0,\zlevel);
\coordinate (C) at (2,2,\zlevel);
\coordinate (D) at (0,2,\zlevel);

% Shadow projected vertically upward and scaled from center (1,1)
\coordinate (As) at ({1 + \scale*(-1)}, {1 + \scale*(-1)}, \shadowz);
\coordinate (Bs) at ({1 + \scale*(1)},  {1 + \scale*(-1)}, \shadowz);
\coordinate (Cs) at ({1 + \scale*(1)},  {1 + \scale*(1)},  \shadowz);
\coordinate (Ds) at ({1 + \scale*(-1)}, {1 + \scale*(1)},  \shadowz);

% Draw the original blue square
\draw[thick,fill=blue!30,opacity=0.8] (A) -- (B) -- (C) -- (D) -- cycle;

% Draw the lighter blue shadow above
\draw[thick,fill=blue!30,opacity=0.8] (As) -- (Bs) -- (Cs) -- (Ds) -- cycle;

% Draw vertical projection lines
\draw[dashed] (A) -- (As);
\draw[dashed] (B) -- (Bs);
\draw[dashed] (C) -- (Cs);
\draw[dashed] (D) -- (Ds);

% Add labels in blue font
\node[blue] at (1,1,\zlevel - 0.3) {$\mathcal{R}(t_1)$};
\node[blue] at (1,1,\shadowz + 0.3) {$\mathcal{R}(t_2)$};
\end{tikzpicture}
\caption{Type (ii) evolution in the absence of reaction and boundary fluxes: the absence of dilution ensures that the \emph{spatial mean} on $\mathcal{R}(t)$ of a spatially uniform species is preserved over time.}
\label{fig:type_2_evolution}
\end{subfigure}
\caption{Pictorial comparison of surface evolution of type (i) and (ii). Color darkness indicates the absolute value of a given species. For type (i) evolution (left), a local increase in surface area generates a dilution effect which lowers concentration locally, in absolute value. For type (ii) evolution (right), there are no dilution effects, as the values of a given species are transported along the material trajectories.}
\end{figure}

\begin{remark}[Interpretation of the modified balance law \eqref{balance_law_eta_type2}]
To give an empirical interpretation of the balance law \eqref{balance_law_eta_type2}, we rewrite it using Reynolds' transport formula:
\begin{equation}
\label{balance_law_eta_type2_no_reaction}
\int_{\mathcal{R}(t)} \material \eta(\boldx,t) \mathrm{d}\boldx =  -\int_{\partial\mathcal{R}(t)} \boldq_\eta \cdot \boldmu + \int_{\mathcal{R}(t)} \rho f(\eta,\theta)\mathrm{d}\boldx, \qquad t\in [0,T].
\end{equation}
In the special case $\boldq_\eta = \0$ and $f=0$,  \eqref{balance_law_eta_type2_no_reaction} becomes
\begin{equation}
\label{balance_law_eta_type2_no_reaction_2}
\int_{\Gamma(t)} \material \eta(\boldx,t) \mathrm{d}\boldx =  0, \qquad t\in [0,T].
\end{equation}
In \eqref{balance_law_eta_type2_no_reaction_2}, if we further assume that $\eta$ is spatially uniform at all times, i.e. $\eta(\boldx,t) = \overline{\eta}(t)$, we obtain
\begin{equation}
|\Gamma(t)| \frac{\mathrm{d}}{\mathrm{d}t}\overline{\eta}(t) = 0, \qquad t\in [0,T],
\end{equation}
which implies $\eta(\boldx,t) = \overline{\eta}(t) = k\in\mathbb{R}$. It follows that
\begin{equation}
\frac{\mathrm{d}}{\mathrm{d}t} \left(\frac{1}{|\Gamma(t)|}\int_{\Gamma(t)} \eta(\boldx,t)\mathrm{d\boldx} \right) = \frac{\mathrm{d}}{\mathrm{d}t} \frac{k|\Gamma(t)|}{|\Gamma(t)|} = 0, \qquad t\in [0,T],
\end{equation}
i.e. the spatial mean of $\eta$ is preserved over time, see Fig. \ref{fig:type_2_evolution} for an illustration. 
\end{remark}

In analogy with the SDIB model \eqref{sdib_model},  in the new model the $\eta$ component will represent the pointwise rate of addition of material \red{onto the surface}; such a physical quantity cannot be viewed as a concentration and is clearly not subject to dilution effects.  Therefore,  $\eta$ will drive type (ii) evolution and will fulfil the balance law \eqref{balance_law_eta_type2}.  By applying integration by parts \eqref{integration_by_parts} to the first term on the right hand side of \eqref{balance_law_eta_type2}, we obtain
\begin{equation}
\label{esdib_eq_eta_derivation_1}
\frac{\mathrm{d}}{\mathrm{d}t} \int_{\mathcal{R}(t)} \eta(\boldx,t)\mathrm{d}\boldx = -\int_{\mathcal{R}(t)} \nablagammaht\cdot\boldq_\eta \mathrm{d}\boldx + \int_{\mathcal{R}(t)} \rho f(\eta,\theta)\mathrm{d}\boldx + \int_{\mathcal{R}(t)} \eta\nabla_{\Gamma(t)} \cdot\dot{\boldx}(t)\mathrm{d}\boldx,
\end{equation}
for $t\in [0,T]$. By using the transport formula \eqref{reynolds} in the left hand side of \eqref{esdib_eq_eta_derivation_1} we obtain
\begin{equation}
\label{esdib_eq_eta_derivation_2}
\int_{\mathcal{R}(t)} \material u_k \mathrm{d}\boldx+ \int_{\mathcal{R}(t)}\nablagammat \cdot \boldq_\eta \mathrm{d}\boldx = \int_{\mathcal{R}(t)} \rho f(\eta,\theta)\mathrm{d}\boldx, \qquad t \in [0,T].
\end{equation}
Since $\mathcal{R}(t)$ is an arbitrary portion of $\Gamma(t)$ we obtain
\begin{equation}
\label{esdib_eq_eta_derivation_3}
\material u_k + \nablagammat \cdot \boldq_\eta = \rho f(\eta,\theta),  \qquad (\boldx,t) \in \mathcal{G}.
\end{equation}
Thanks to Fickian diffusion \eqref{fick} (motivated as in the SDIB model \eqref{sdib_model}, see \cite{lacitignola2017turing}), \eqref{esdib_eq_eta_derivation_3} finally yields the equation for $\eta$:
\begin{equation}
\label{esdib_eq_eta_derivation_4}
\partial^\bullet \eta - \Delta_{\Gamma(t)} \eta = \rho f(\eta,\theta), \qquad (\boldx,t) \in \mathcal{G}.
\end{equation}
The $\theta$ component is the \red{surface} concentration of a surfactant confined to the surface $\Gamma(t)$ at all times; $\theta$ will therefore drive type (i) evolution and will fulfil the standard balance law \eqref{balance_law_theta}.  Therefore, the equation for $\theta$ will coincide with the second equation of the ESRDS \eqref{esrds}:
\begin{equation}
\label{esdib_eq_theta_derivation}
\partial^\bullet \theta + \theta \nabla_{\Gamma(t)} \cdot \dot{\boldx}(t) - d\Delta_{\Gamma(t)} \theta = \rho g(\eta, \theta),\qquad (\boldx,t) \in\mathcal{G}.
\end{equation}
By combining \eqref{evolution_law}, \eqref{esdib_eq_eta_derivation_4}, and \eqref{esdib_eq_theta_derivation}, we finally obtain the non-dimensionalised Evolving Surface DIB (ESDIB) model:
\begin{equation}
\label{model_evolving}
\begin{cases}
\partial^\bullet \eta - \Delta_{\Gamma(t)} \eta = \rho f(\eta,\theta), \qquad (\boldx,t) \in\mathcal{G};\\
\partial^\bullet \theta + \theta \nabla_{\Gamma(t)} \cdot \dot{\boldx}(t) - d\Delta_{\Gamma(t)} \theta = \rho g(\eta, \theta),\qquad (\boldx,t) \in\mathcal{G};\\
\red{\nablagammat\eta \cdot \boldmu = 0, \qquad (\boldx,t) \in\mathcal{G}_\partial;}\\
\red{\nablagammat\theta \cdot \boldmu = 0, \qquad (\boldx,t) \in\mathcal{G}_\partial;}\\
\dot{\boldx} = \kappa \eta \boldn, \qquad (\boldx,t) \in\mathcal{G};\\
\eta(\cdot,0) = \eta_0,  \qquad \boldx \in \Gamma(0);\\
\theta(\cdot,0) = \theta_0, \qquad \boldx \in \Gamma(0);\\
\Gamma(0) = \Gamma_0.
\end{cases}
\end{equation}
\red{We remark that the derivation of the above model \eqref{model_evolving} holds true regardless of the functional form of the kinetics $(f,g)$, as long as the balance laws \eqref{balance_law_theta} and \eqref{balance_law_eta_type2} are fulfilled.  
The ESDIB model \eqref{model_evolving} is presented in non-dimensional form. The reader interested in its dimensional counterpart is referred to our previous work \cite[Section 8]{lacitignola2017turing} on the SDIB model \eqref{sdib_model}.  The adimensionalisation presented therein holds true for the ESDIB model \eqref{model_evolving}, as well.   We also remark that, depending on the initial domain, initial condition and model parameters, there might exist a time $t_{\text{self}} > 0$ where a self-intersection occurs in the surface $\Gamma(t)$, as illustrated in Fig. \ref{fig:esdib}.  For $t\geq t_{\text{self}}$, the ESDIB model \eqref{model_evolving} is still well-posed mathematically but not physically.  In the experimental setting, self-intersections correspond to metallic protrusions merging with each other,  leading to the formation of sponge-like structures.  The mathematical translation of this phenomenon is a topological change of $\Gamma(t)$ for $t = t_{\text{self}}$.  In the ESDIB model \eqref{model_evolving}, instead, the surface $\Gamma(t)$ is intangible and remains topologically equivalent to the initial surface $\Gamma_0$ even after self-touching and self-crossing, see Fig. \ref{fig:esdib}.  In Section \ref{sec:numerical_examples} we notice, however, that the ESDIB model shows good qualitative agreement with the experiments even after self-intersections.  The incorporation of topological changes in the ESDIB model is an open research direction. }

\section{Spatial discretisation}
\label{sec:space_discretisation}
In this Section, we derive the weak formulation of the ESDIB model \eqref{model_evolving} and we present its spatial discretisation using the Evolving Lumped Surface Finite Element Method (LESFEM) \cite{frittelli2018numerical}.

\subsection{Weak formulation}
By multiplying the first two equations in \eqref{model_evolving} by two test functions $\phi,\psi$ and integrating over $\Gamma(t)$ we obtain
\begin{equation}
\label{weak_1}
\begin{cases}
\displaystyle\int_{\Gamma(t)} \varphi \material \eta - \intgammat \varphi\deltagammat \eta = \rho\intgammat f(\eta,\theta)\varphi;\\
\displaystyle\int_{\Gamma(t)} \psi \material \theta + \int_{\Gamma(t)}\theta\psi \nabla_{\Gamma(t)} \cdot \dot{\boldx}(t) - d\intgammat \psi\deltagammat \theta \!=\! \rho\intgammat g(\eta,\theta)\psi.
\end{cases}
\end{equation}
By applying Green's formula on surfaces \cite{dziuk2013finite} in both equations in \eqref{weak_1} and Reynolds' transport \red{formula \eqref{reynolds}} in the second equation of \eqref{weak_1},  system \eqref{weak_1} becomes
\begin{equation}
\label{weak_2}
\begin{cases}
\displaystyle\int_{\Gamma(t)} \varphi \material \eta + \intgammat \nablagammat \varphi\cdot\nablagammat \eta = \rho \intgammat f(\eta,\theta)\varphi;\\
\displaystyle\frac{\mathrm{d}}{\mathrm{d}t}\left(\intgammat\varphi\theta\right) \!-\!\intgammat  \theta\material \varphi \!+\! d\intgammat \nablagammat\psi\cdot \nablagammat \theta \!=\! \rho\intgammat \!g(\eta,\theta)\psi.
\end{cases}
\end{equation}

\begin{remark}
The weak formulation of \red{the} usual evolving \red{surface} RDSs of the form \eqref{esrds}, where surface evolution is of type (i) in both equations, is (see \cite{frittelli2018numerical}):
\begin{equation}
\begin{cases}
\vspace{1mm}
\displaystyle \frac{\mathrm{d}}{\mathrm{d}t}\left(\intgammat\varphi\eta\right) \!-\!\intgammat  \eta\material \varphi  \!+\! \intgammat \nablagammat\varphi\cdot\nablagammat\eta \!=\! \rho\intgammat f(\eta,\theta)\varphi;\\
\displaystyle \frac{\mathrm{d}}{\mathrm{d}t}\left(\intgammat\varphi\theta\right) \!-\!\intgammat  \theta\material \varphi  \!+\! d\intgammat \nablagammat\psi\cdot\nablagammat\theta \!=\! \rho\intgammat \!g(\eta,\theta)\psi.
\end{cases}
\end{equation}
\end{remark}

\subsection{Lumped Evolving Surface Finite Elements}
We introduce some notation following \cite{frittelli2018numerical}. Let $\Gamma_h(t)$ be a time-dependent triangulation of the evolving surface $\Gamma(t)$. Let $\mathbb{V}_h(t)$ be the space of piecewise linear functions on $\Gamma_h(t)$ and let $\{\chi_i(\cdot,t)\}_{i=1}^N$ be the Lagrangian basis of $\mathbb{V}_h(t)$.  Let $I_h:\mathcal{C}^0(\Gamma_h(t)) \rightarrow \mathbb{V}_h(t)$ be the piecewise linear interpolant operator.  The discrete material derivative is defined by
 \begin{equation}
\label{discrete_material_derivative}
\materialh U(\boldx,t) := \frac{\partial \widetilde{U}}{\partial t} + I_h(\boldv)\cdot \nabla \widetilde{U},
\end{equation}
where $\widetilde{U}: \mathcal{N}_h(t)\rightarrow\mathbb{R}$ is a continuous and piecewise differentiable extension of $U$ defined on an open neighborhood $\mathcal{N}_h(t)$ of $\Gamma_h(t)$.
The spatially discrete formulation of \eqref{weak_2} is given by
\begin{equation}
\label{spatially_discrete}
\begin{cases}
\displaystyle\intgammaht I_h\left(\chi_i \materialh \eta_h\right) + \intgammaht \nablagammaht \eta_h\cdot\nablagammaht  \chi_i = \rho\intgammaht I_h(f(\eta_h,\theta_h) \chi_i);\\
\displaystyle\frac{\mathrm{d}}{\mathrm{d}t}\intgammaht \!\!\!I_h\left(\theta_h  \chi_i\right) \!-\!\intgammaht \theta_h\material_h \chi_i \!+\! d\intgammaht \nablagammaht\theta_h\cdot \nablagammaht  \chi_i \!=\! \rho\intgammaht \!\!\! I_h\left(g(\eta_h,\theta_h) \chi_i\right),
\end{cases}
\end{equation}
for all $i=1,\dots,N$. We express $\eta_h$ and $\theta_h$ in the Lagrange basis as follows:
\begin{equation}
\label{eta_h_theta_h}
\eta_h(\boldx,t) = \sum_{j=1}^N a_j(t) \chi_j(\boldx,t), \qquad \theta_h(\boldx,t) = \sum_{j=1}^N b_j(t) \chi_j(\boldx,t),  \qquad t\in [0,T],  \ \boldx \in \Gamma_h(t).  
\end{equation}
\red{We recall the following fundamental property of the Lagrangian basis functions.
\begin{lemma}[Transport property of the Lagrangian basis functions, \cite{dziuk2007finite}]
The Lagrangian basis functions $\chi_i$, $i=1,\dots,N$, fulfil
\begin{equation}
\label{transport_property_basis}
\material_h \chi_i = 0, \qquad i=1,\dots,N. 
\end{equation}
Moreover,  the material derivative of $\eta_h$ and $\theta_h$ in \eqref{eta_h_theta_h} fulfils
\begin{equation}
\label{material_eta_h_theta_h}
\materialh\eta_h(\boldx,t) = \sum_{j=1}^N \dot{a}_j(t) \chi_j(\boldx,t), \qquad \materialh\theta_h(\boldx,t) = \sum_{j=1}^N \dot{b}_j(t) \chi_j(\boldx,t),  \qquad t\in [0,T], \ \boldx \in \Gamma_h(t).
\end{equation}
\end{lemma}}
Using \eqref{transport_property_basis}-\eqref{material_eta_h_theta_h},  the spatially discrete formulation \eqref{spatially_discrete} takes the following matrix-vector form:
\begin{equation}
\label{spatially_discrete_1}
\begin{cases}
\displaystyle\sum_{j=1}^N \dot{a}_j\intgammaht I_h\left(\chi_j \chi_i\right) + \sum_{j=1}^N a_j\intgammaht \nablagammaht \chi_j\cdot\nablagammaht  \chi_i = \rho f(a_j,b_j)\intgammaht I_h(\chi_j\chi_i);\\
\displaystyle\frac{\mathrm{d}}{\mathrm{d}t}\!\left(\sum_{j=1}^N a_j\!\intgammaht \!\!\!\!I_h\left(\chi_j  \chi_i\right)\right) \!+\! d\sum_{j=1}^N a_j\! \intgammaht\!\!\!\! \nablagammaht\chi_j\cdot \nablagammaht  \chi_i \!=\! \rho g(a_j,b_j)\!\intgammaht \!\!\!\!I_h\left(\chi_j \chi_i\right),
\end{cases}
\end{equation}
for all $i=1,\dots,N$.  Using the lumped mass matrix $M\in\mathbb{R}^{N\times N}$ defined by
\begin{equation}
M_{ij} = I_h(\chi_i\chi_j) = \begin{cases}
\displaystyle\int_{\Gamma(t)} \chi_i \qquad &\text{if } i=j;\\
0 & \text{if } i \neq j,
\end{cases}
\end{equation}
the system \eqref{spatially_discrete_1} takes the matrix-vector form
\begin{equation}
\label{matrix-vector-form}
\begin{cases}
\vspace{1mm}
\dfrac{\mathrm{d}}{\mathrm{d}t}\boldx_i(t) = \kappa a_i(t) \boldn_i(t), \qquad i=1,\dots,N;\\
\vspace{1mm}
M(t)\dfrac{\mathrm{d}}{\mathrm{d}t}\bolda(t) + K(t) \bolda(t) = \rho M(t) f(\bolda(t),\boldb(t));\\
\dfrac{\mathrm{d}}{\mathrm{d}t}\left(M(t)\boldb(t)\right) + K(t) \boldb(t) = \rho M(t) g(\bolda(t),\boldb(t)),
\end{cases}
\qquad t\in [0,T], 
\end{equation}
where $\boldn_i(t)$ is the outward unit normal vector on $\Gamma_h(t)$ evaluated at the node $\boldx_i(t)$.

\begin{remark}[Discretisation of usual evolving RDSs]
The spatial discretisation of usual evolving RDSs of the form \eqref{esrds} is (see \cite{frittelli2018numerical}):
\begin{equation}
\begin{cases}
\vspace{1mm}
\dfrac{\mathrm{d}}{\mathrm{d}t}\boldx_i(t) = \kappa a_i(t) \boldn_i(t), \qquad i=1,\dots,N;\\
\vspace{1mm}
\dfrac{\mathrm{d}}{\mathrm{d}t}\left(M(t)\bolda(t)\right) + K(t) \bolda(t) = \rho M(t) f(\bolda(t),\boldb(t));\\
\dfrac{\mathrm{d}}{\mathrm{d}t}\left(M(t)\boldb(t)\right) + K(t) \boldb(t) = \rho M(t) g(\bolda(t),\boldb(t)),
\end{cases}
\qquad t\in [0,T].
\end{equation}
\end{remark}

\section{Time discretisation}
\label{sec:time_discretisation}
We choose a discretisation step $h_t > 0$,  we define $N_T := \lceil \frac{T}{h_t}\rceil$ and, for all $n=0,\dots,h_t$ we define the $n$-th time step as $t_n := nh_t$.
We adopt the following time discretisation of \eqref{matrix-vector-form}
\begin{equation}
\label{imex-euler}
\begin{cases}
\dfrac{\boldx_{i,n+1} - \boldx_{i,n}}{\tau} = \boldv_{i,n} = \kappa a_{i,n} \boldn_{i,n}, \qquad i=1,\dots,N;\\
M_n\dfrac{\bolda_{n+1}-\bolda_n}{\tau} + K_{n+1} \bolda_{n+1} =  \rho M_n f(\bolda_n,\boldb_n);\\
\dfrac{M_{n+1}\boldb_{n+1}-M_n \boldb_n}{\tau} + K_{n+1} \boldb_{n+1} = \rho M_n g(\bolda_n,\boldb_n).
\end{cases}
\end{equation}
In \eqref{imex-euler}, diffusion is approximated implicitly (hence the stiffness matrix $K_{n+1}$ evaluated at $t_{n+1}$), while reactions are evaluated implicitly (hence the mass matrix $M_{n+1}$ evaluated at $t_{n+1}$ on the right-hand side of the second and third equations).  This choice is called the Implicit-Explicit (IMEX) Euler discretisation, and was successfully applied for the time discretisation of surface PDEs in previous works \cite{lacitignola2017turing,lacitignola2019spiral,frittelli2018numerical,frittelli2017lumped, frittelli2019preserving} since:
\begin{itemize}
\item it is unconditionally stable for the (stiff) diffusion part when the surface is stationary \cite{frittelli2017lumped, frittelli2019preserving} or evolves according to known isotropic growth \cite{frittelli2018numerical};
\item it is implementation-firiendly, since the second and third equations in \eqref{imex-euler} are linear algebraic systems to be solved at each time step, thereby avoiding the nonlinear rootfinding step usually required by implicit time discretisations.
\end{itemize} 
In the evolving surface framework, however, node position needs to be computed at each time step by \red{solving for the discrete surface}. In the scheme \eqref{imex-euler}, the first equation is obtained through Explicit Euler discretisation of \eqref{evolution_law}, in order to avoid a nonlinear rootfinding step (since each $\boldn_{i,n+1}$ is a function of $\boldx_{i,n+1}$ itself and of its neighboring nodes).  The matrix $M_{n+1}$ appearing on the left hand side of the third equation of \eqref{imex-euler} is then computed using the $\boldx_{i,n+1}$'s.  This particular adaptation of the IMEX Euler method was employed in \cite{kovacs2017convergence},  where it was referred to as ``a linearly implicit Euler'' scheme.  In that work, the surface undergoes type (i) growth for both species.  Here instead, the second equation in \eqref{imex-euler} is different due to type (ii) growth for the species $\eta$.  

\section{Numerical examples and experimental comparisons}
\label{sec:numerical_examples}
In this section we 
\begin{itemize}
\item showcase numerical simulations of the ESDIB model for different choices of the parameters and initial domain
\item compare the numerical solutions of the ESDIB and DIB models on equal parameters and spatial domain (initial spatial domain for the ESDIB model)
\item compare the numerical solutions of the ESDIB model with experimental samples.
\end{itemize}
\red{Following our previous work \cite{sgura2019parameter} on the DIB model,  all the model values except $B$, $C$, and $D$ are chosen as follows in all the numerical examples:
\begin{equation}
d = 20, \ \alpha = 0.5, \ \gamma = 0.2, \ k_2 = 2.5, \ k_3 = 1.5, \ A_1 = 10, \ A_2 = 1, \ k=0.2.
\end{equation}
The values of $B$ and $C$ are specified in each experiment and are summarized in Table \ref{tab:summary_experiments}.  Finally, the value of $D$ is determined each time as $D = C\frac{(1-\alpha)(1-\gamma+\gamma\alpha)}{\alpha(1+\gamma\alpha)}$.}
\begin{table}[ht!]
\centering
\red{\begin{tabular}{c|l|c|c|c}
Example no.  & Initial domain $\Gamma_0$ & $B$ & $C$ & DIB pattern class\\
\hline\hline
1 & Square with edge $L=20$ & 30 & 3 & holes\\
2 & Square with edge $L=30$ & 66 & 3 & labyrinth\\
3 & Square with edge $L=20$ & 62 & 5 & spots\\
4 & Sphere with radius $R = 3$ & 30 & 3 & holes\\
5 & Sphere with radius $R = 10$ & 66 & 3 & labyrinth\\
6 & Sphere with radius $R = 5$ & 62 & 5 & spots
\end{tabular}}
\caption{\red{Summary of the numerical experiments.}}
\label{tab:summary_experiments}
\end{table}

\noindent
\red{In all the experiments, the initial condition is the following random perturbation of the equilibrium \eqref{equilibrium_dib}:
\begin{equation}
\eta_0(\boldx) = \eta_{\text{eq}} + 10^{-4}r(\boldx), \qquad \theta_0(\boldx) = \theta_{\text{eq}} + 10^{-4}r(\boldx),
\end{equation}
for all $\boldx\in\Omega$ in the DIB model \eqref{dib_model} or $\boldx\in\Gamma_0$ in the ESDIB model \eqref{model_evolving}, where $r(\boldx)$ is a spatial random perturbation with uniform distribution over $[-1,1]$ and $\eta_{\text{eq}}$ and $\theta_{\text{eq}}$ are the uniform steady states defined in \eqref{equilibrium_dib}.  We remark that,  both in the DIB \eqref{dib_model} and the ESDIB \eqref{model_evolving} models,  negative values for $\eta$ are physically meaningful and correspond to areas where corrosion is taking place.}
\red{Since, in the ESDIB model, surface evolution causes the mesh to become increasingly irregular over time,  the stopping time is determined in each experiment in such a way to avoid the formation of degenerate elements.  To prevent this issue,  an ALE evolving mesh could be employed \cite{elliott2012ale}. In this work, however, we were able to obtain satisfactory experimental comparisons without relying on remeshing techniques, in order to prioritise computational efficiency on the limited available computational architecture. The timestep is $\tau = 1e$-2 in all the experiments. The computations were carried out in MATLAB R2024b.}

\subsection{Experimental cases}
The experiments we are considering here concern the morphochemical evolution of two kinds of zinc-battery anodes --zinc foils and zinc sponges--, operated in different electrochemical regimes and analyzed post mortem at the mesoscale by scanning electron microscopy (SEM). The choice of these anodes types resides in the fact that they are currently the two key options for the research and development of next-generation, post-lithium batteries of the sealed type. From the mathematical side, these two options dictate the choice of these \red{geometries} of the initial domain: planar for the zinc foil case and spherical for the zinc-sponge one: the latter choice is due to the fact that the zinc units resulting from the material synthesis are micrometric spheroids, \red{unless otherwise specified}. In this research, zinc foil anodes are operated in symmetric cells with an alkaline aqueous electrolyte without and with additives: details on electrochemical materials science and cycling are available in \cite{rossi2020electrodeposition}.  As far as zinc-sponge anodes are concerned, full details on material preparation, electrochemical testing and 3D imaging, both ex situ and in operando, as well as a comprehensive account of the literature in the field, can be found in \cite{bozzini2020morphological, bozzini2025formation}. The formulation of zinc-sponge anodes is one of the currently considered microstructure tailoring approaches, aimed at achieving extensive charge-discharge cycling for batteries with metal anodes. A zinc sponge consists of a framework of connected metallic branches, coated with a layer of ZnO. In their anodic operation, continuity of the electron-conductive metallic network is presupposed to persists down to the required depth of discharge and to be preserved over cycling. Moreover,  Zn$^{2+}$ is confined into the porous electrode structure, preventing unstable metal growth. Microscopic imaging of the actual distribution of zinc and ZnO, resulting from different operating conditions and the capability of modelling it is crucial for the understanding of the actual behaviour of these materials and for their knowledge-based improvement. In the present study, we shall report original experimental work referring to the material, synthesized according to the protocols of \cite{rossi2020electrodeposition, bozzini2020morphological} and charged/cycled in realistic battery conditions, that will be detailed where relevant. In most cases (with a few exception indicated below), the last electrochemical step is recharge, corresponding to electrochemical growth of the metal phase.

\subsection{Example 1: square, $B=30$, $C=3$, holes}
\begin{figure}[ht!]
\begin{subfigure}{\textwidth}
\vspace*{-5mm}
\includegraphics[scale=0.35]{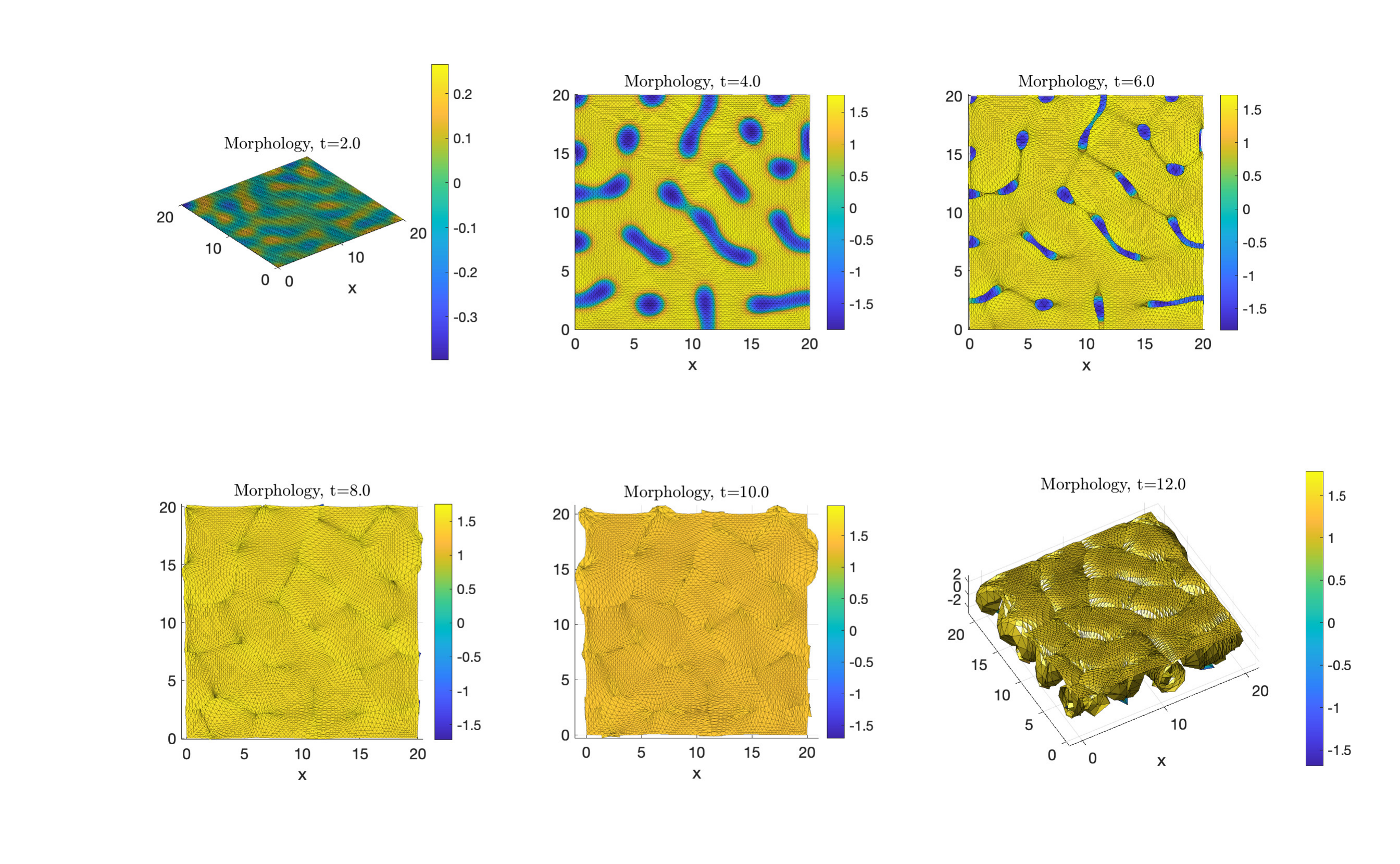}
\caption{}
\label{fig:example_1_esdib_num}
\end{subfigure}
\begin{subfigure}{\textwidth}
\centering
\includegraphics[scale=.25]{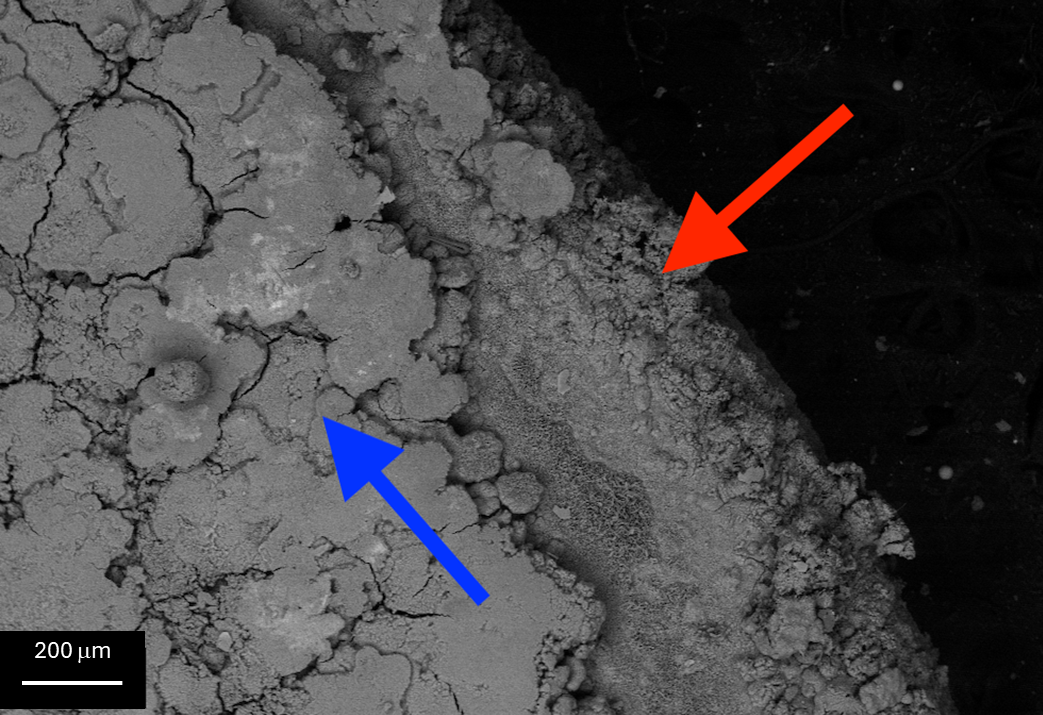}
\caption{}
\label{fig:example_1_esdib_lab}
\end{subfigure}
\caption{Example 1.  (a) $\eta$ component of the numerical solution of the ESDIB model at various times.  (b) SEM micrograph of a zinc foil electrode charged in pure 6M KOH solution at 5 mA cm$^{-2}$ for 3 hours.  Blue and red arrows indicate low and high current density regions, respectively.}
\label{fig:example_1_esdib}
\end{figure}

\begin{figure}[ht!]
\hspace*{-5mm}
\includegraphics[scale=0.25]{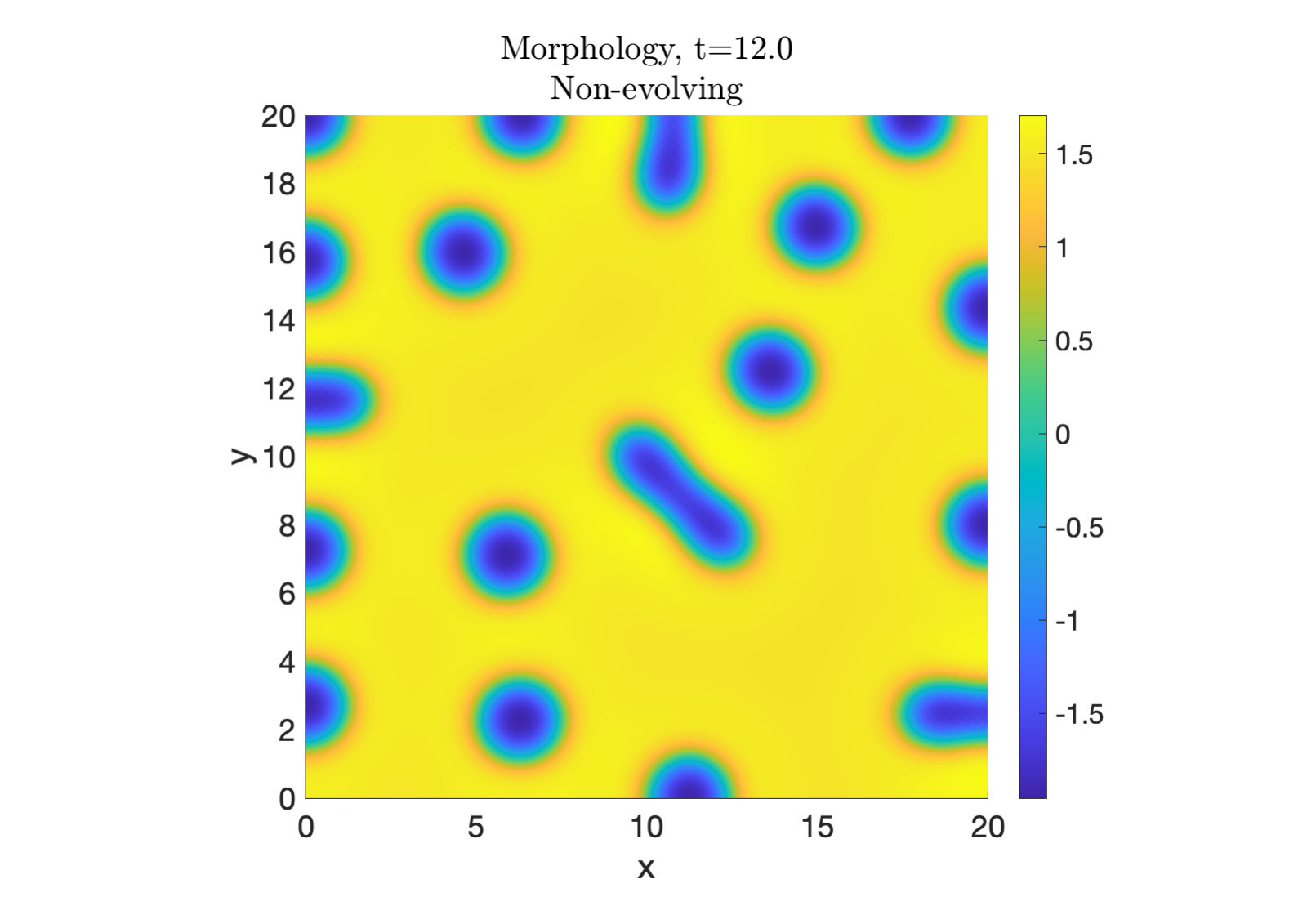}\hspace*{-10mm}
\includegraphics[scale=0.3]{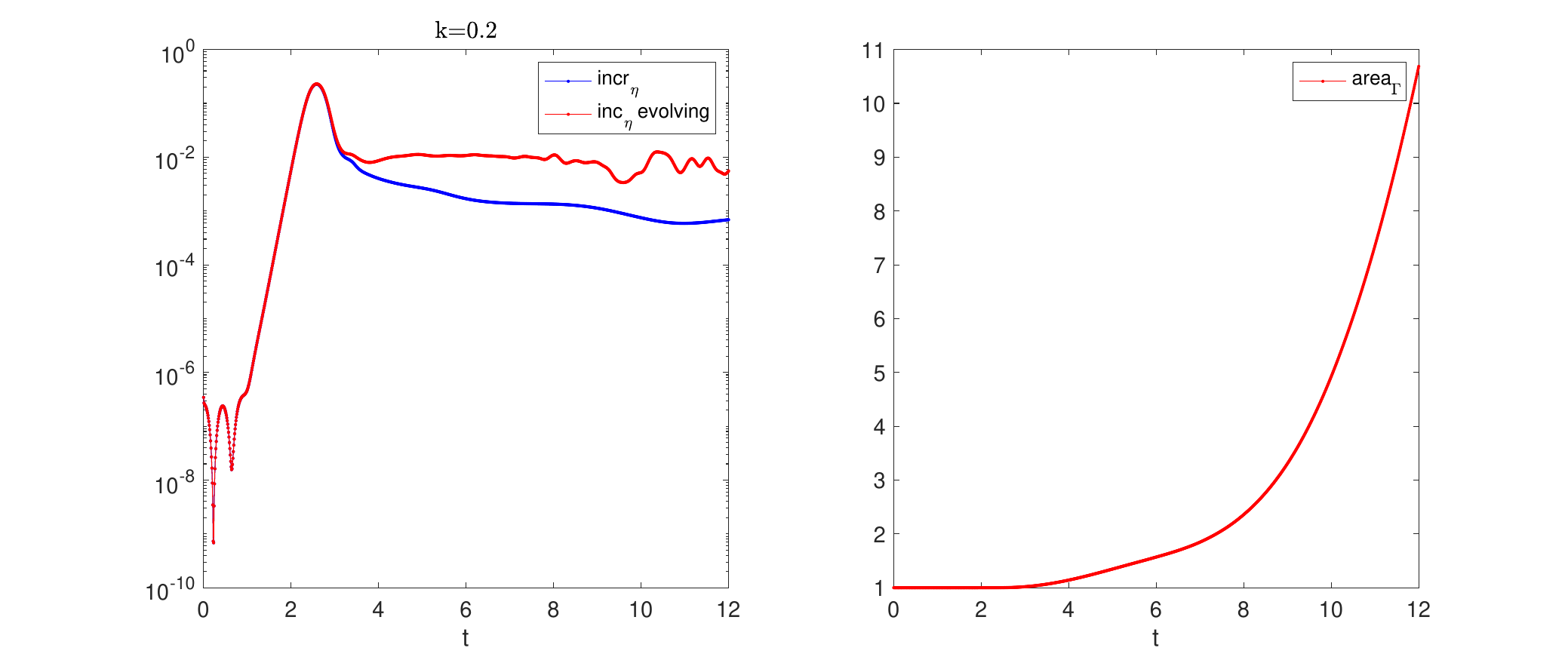}
\caption{Example 1. Left: \red{$\eta$ component} of the DIB model at the final time.  Middle: time increment \red{$\|\eta(t_{i+1}-t_i)\|_2$ of the $\eta$ component} of the DIB and ESDIB models.  Right: area of the evolving surface $\Gamma(t)$ in the ESDIB model.}
\label{fig:example_1_dib_and_increment}
\end{figure}

In this example, the initial domain is the square $\Omega = [0,L]^2$ for $L=20$. The final time is $T=12$. The \red{model}  parameters are chosen as $B = 30$, $C=3$,  \red{which lead to a reversed spots (holes)-type pattern in the DIB model \cite{sgura2019parameter}}, and $\rho=2$.  In Figure \ref{fig:example_1_esdib}, we show the numerical solution of the ESDIB model at different times (subfigure \ref{fig:example_1_esdib_num}) and we compare it with the SEM micrograph of a zinc foil electrode charged in pure 6M KOH solution at 5 mA cm$^{-2}$ for 3 hours (subfigure \ref{fig:example_1_esdib_lab}).  Specifically,  (subfigure \ref{fig:example_1_esdib_lab}) we show the edge of the electrode: flat crystallites ca. 10$\mu$m in size grow in the internal, lower current-density region of the electrode (indicated by a red arrow), while smaller grains are found at the electrode edge, where a higher current density develops (indicated by a blue arrow). In fact, a higher current density, correlates with a higher nucleation rate, ultimately leading to a finer texture \cite{bozzini2012morphogenesis}. In Figure \ref{fig:example_1_dib_and_increment},  we show the solution of the DIB model at the final time, the time increments of the solutions of both models and the area of $\Gamma(t)$ in the ESDIB model over time.  \red{The exponential growth of such area reflects the early stages of structure development and doubling/branching.  Physically, such growth ceases to be exponential when self-intersections and subsequent topological changes begin to occur.  As discussed previously, such topological changes are not captured by the ESDIB model and therefore we observe exponential surface growth at all times. This holds true in all the remaining experiments. Finally,  the $\theta$ component exhibits a similar behaviour as the $\eta$ component and is therefore not omitted from the plots in all the experiments. }

\subsection{Example 2: square, $B=66$, $C=3$, labyrinth}
\begin{figure}[p]
\begin{center}
\begin{subfigure}{\textwidth}
\vspace*{-10mm}
\includegraphics[scale=0.35]{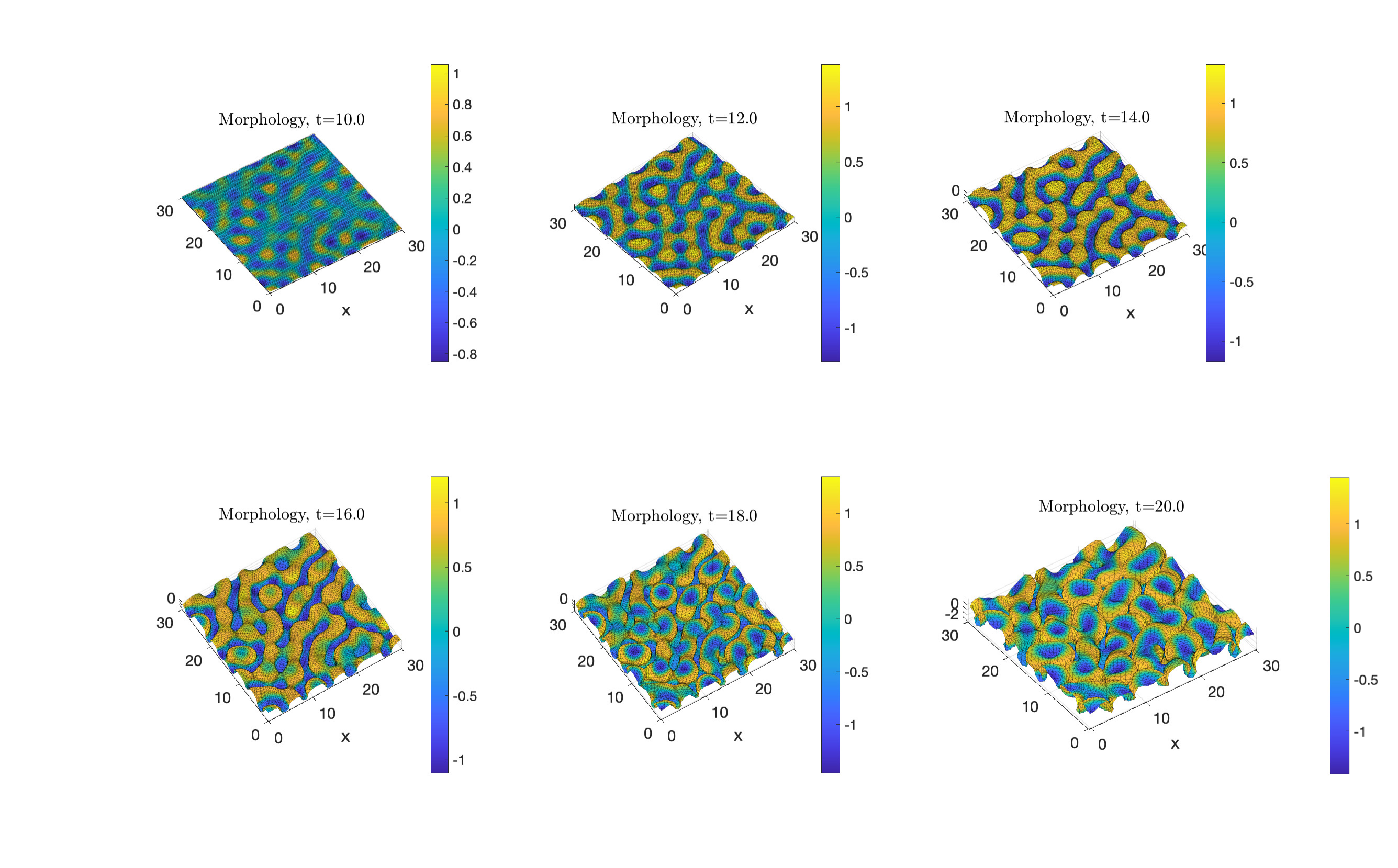}
\caption{}
\label{fig:example_2_esdib_num}
\end{subfigure}
\begin{subfigure}{\textwidth}
\centering
\includegraphics[scale=.25]{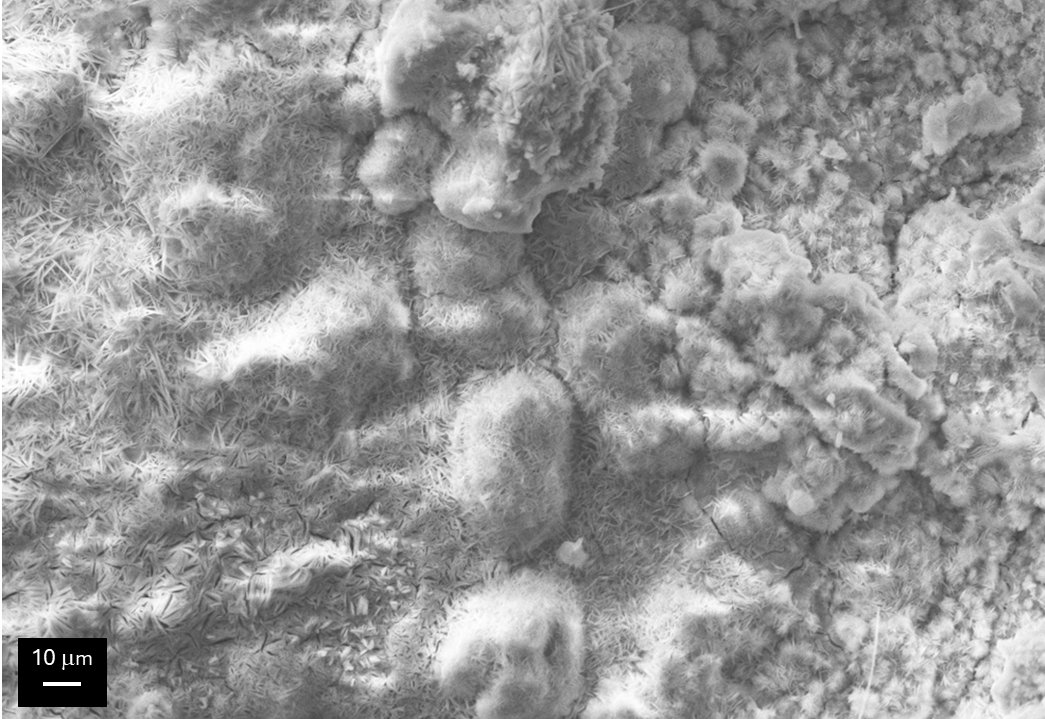}
\caption{}
\label{fig:example_2_esdib_lab}
\end{subfigure}
\end{center}
\caption{Example 2.  (a) $\eta$ component of the numerical solution of the ESDIB model at various times.  (b) SEM micrograph of a zinc foil electrode charged in a 6M KOH solution containing 100 ppm cetyl-trimethyl ammonium chloride additive, at 20 mA cm$^{-2}$ for 30 min.}
\label{fig:example_2_esdib}
\end{figure}

\begin{figure}[p]
\centering
\includegraphics[scale=0.3]{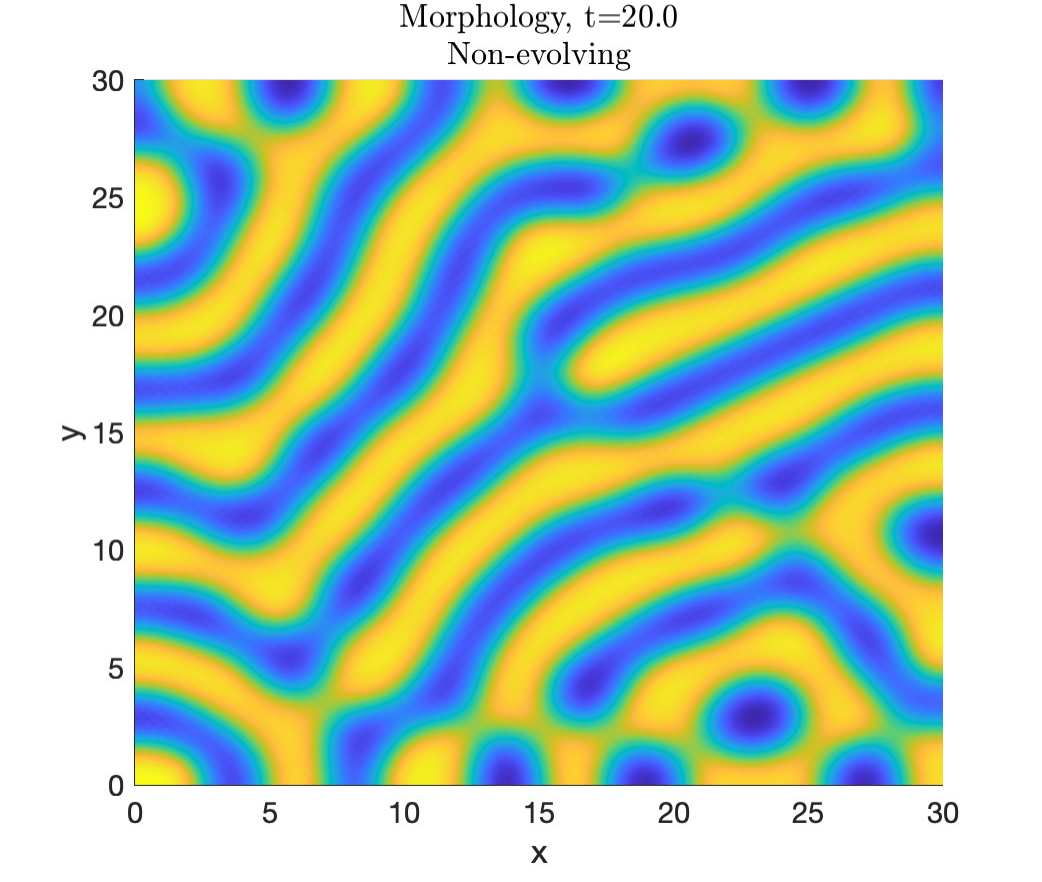}\hspace*{-10mm}
\includegraphics[scale=0.28]{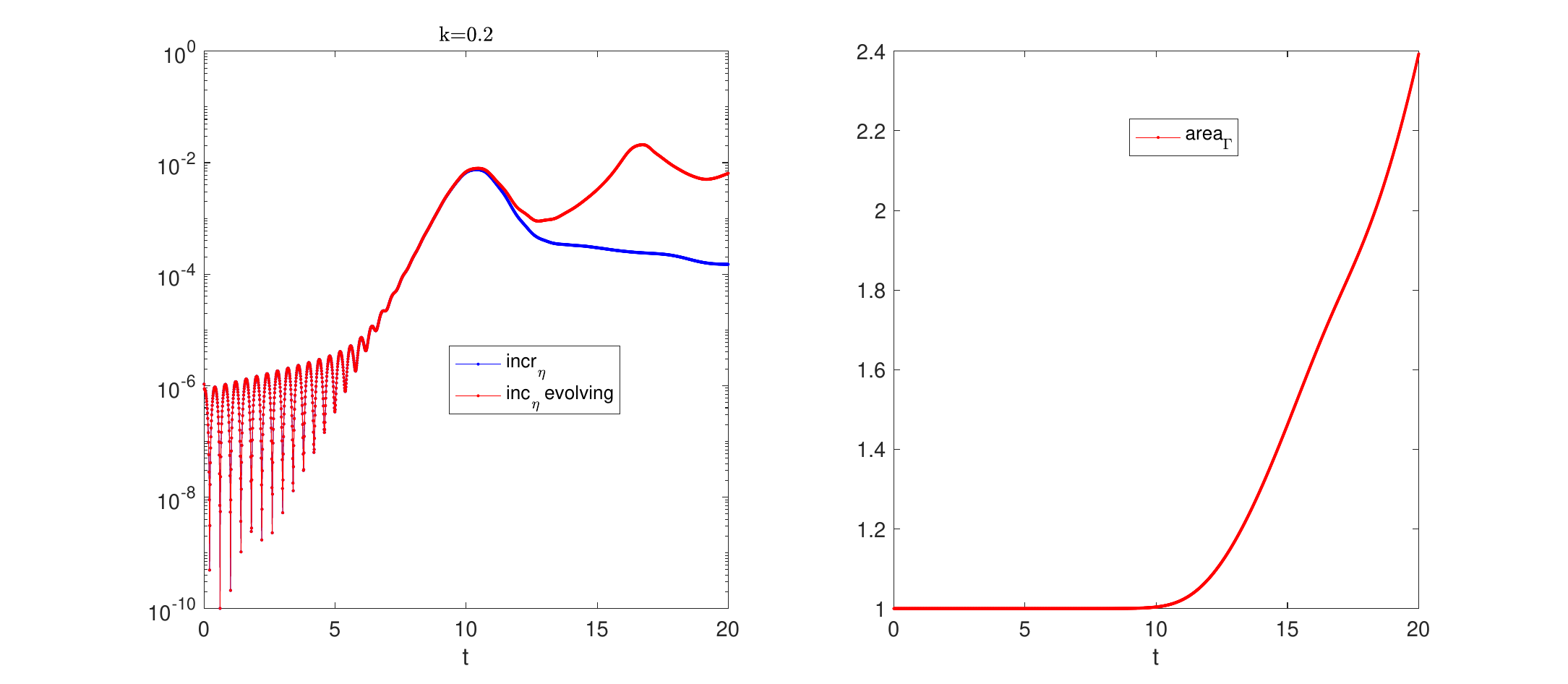}
\caption{Example 2. Left: $\eta$ component of the DIB model at the final time.  Middle: time increment of the numerical solutions of the DIB and ESDIB models.  Right: area of the evolving surface $\Gamma(t)$ in the ESDIB model.}
\label{fig:example_2_dib_and_increment}
\end{figure}

In this example, the initial domain is the square $\Omega = [0,L]^2$ for $L=30$. The final time is $T=20$. The \red{model} parameters are chosen as $B = 66$, $C=3$, \red{which lead to a labyrinth-type pattern in the DIB model \cite{sgura2019parameter}}, and $\rho=1$.  In Figure \ref{fig:example_2_esdib}, we show the numerical solution of the ESDIB model at different times (subfigure \ref{fig:example_2_esdib_num}) and we compare the numerical solution with the SEM morphology of a zinc foil electrode charged in a 6M KOH solution containing 100 ppm cetyl-trimethyl ammonium chloride additive, at 20 mA cm$^{-2}$ for 30 min (subfigure \ref{fig:example_2_esdib_lab}). This additive favours the formation of a partial surface coverage with insulated compact zinc grains of globular type, separared by flat nucleation-exclusion zones \cite{bozzini2012morphogenesis} of the dimensions comparable to those of the grains \cite{sgura2017xrf}.  In Figure \ref{fig:example_2_dib_and_increment},  we show the solution of the DIB model at the final time, the time increments of the solutions of both models and the area of $\Gamma(t)$ in the ESDIB model over time.

\subsection{Example 3: square, $B=62$, $C=5$, spots}
\begin{figure}[p]
\begin{subfigure}{\textwidth}
\includegraphics[scale=0.35]{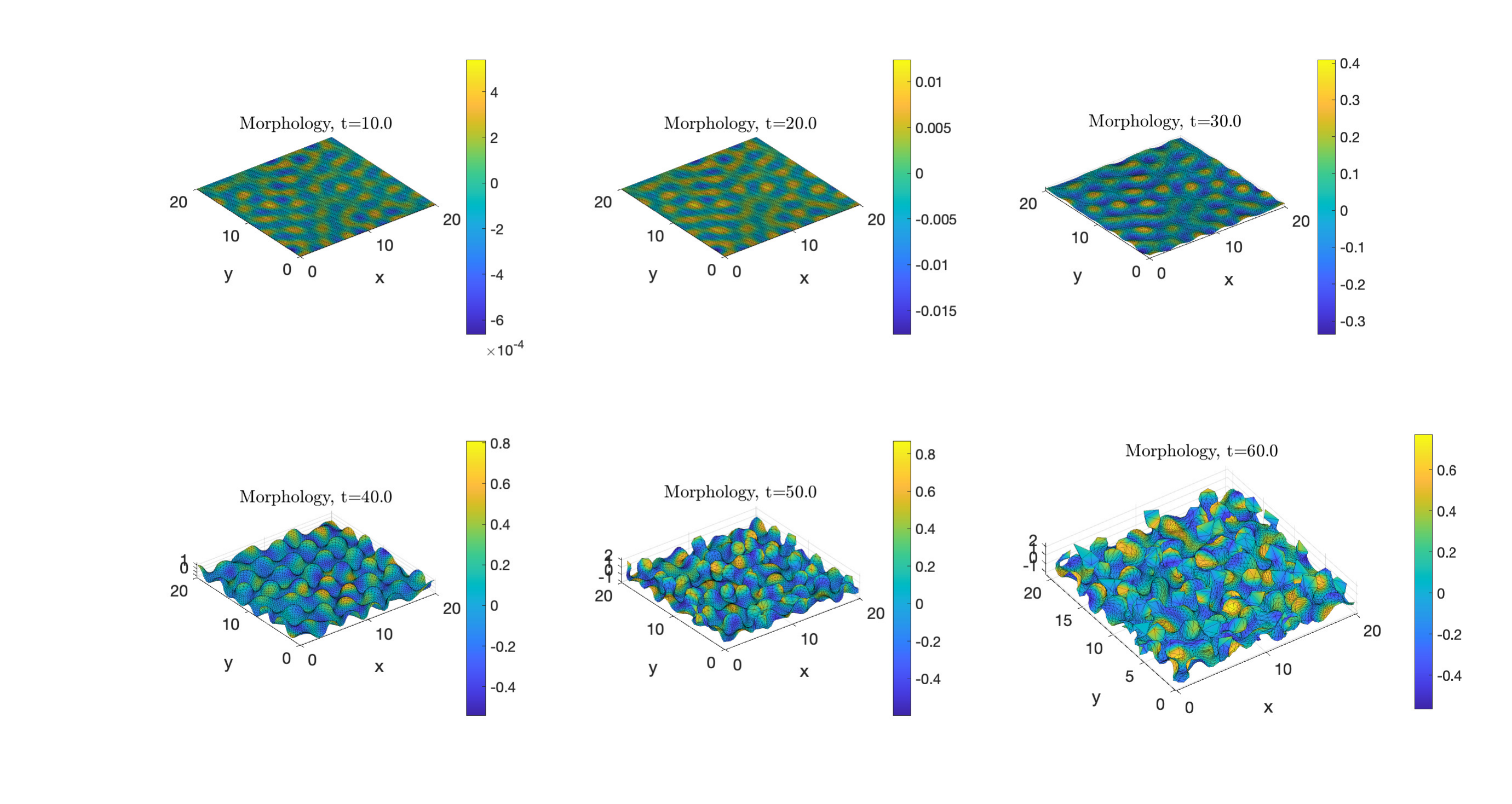}
\caption{}
\label{fig:example_3_esdib_num}
\end{subfigure}
\begin{subfigure}{\textwidth}
\centering
\includegraphics[scale=.25]{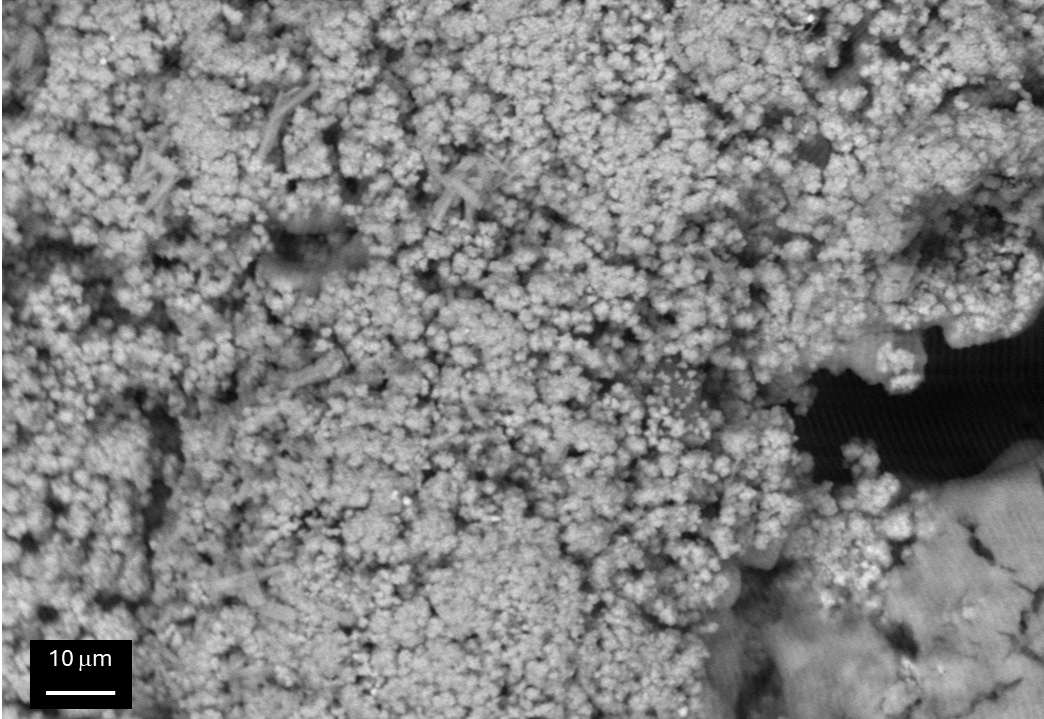}
\caption{}
\label{fig:example_3_esdib_lab}
\end{subfigure}
\caption{Example 3. (a) $\eta$ component of the numerical solution of the ESDIB model at various times. (b) SEM micrograph of a zinc foil anode charged in a 6M KOH solution containing 100 ppm of tetra-butyl ammonium bromide at 5 mA cm$^{-2}$ for 1 hour.}
\label{fig:example_3_esdib}
\end{figure}

\begin{figure}[p]
\centering
\includegraphics[scale=0.27]{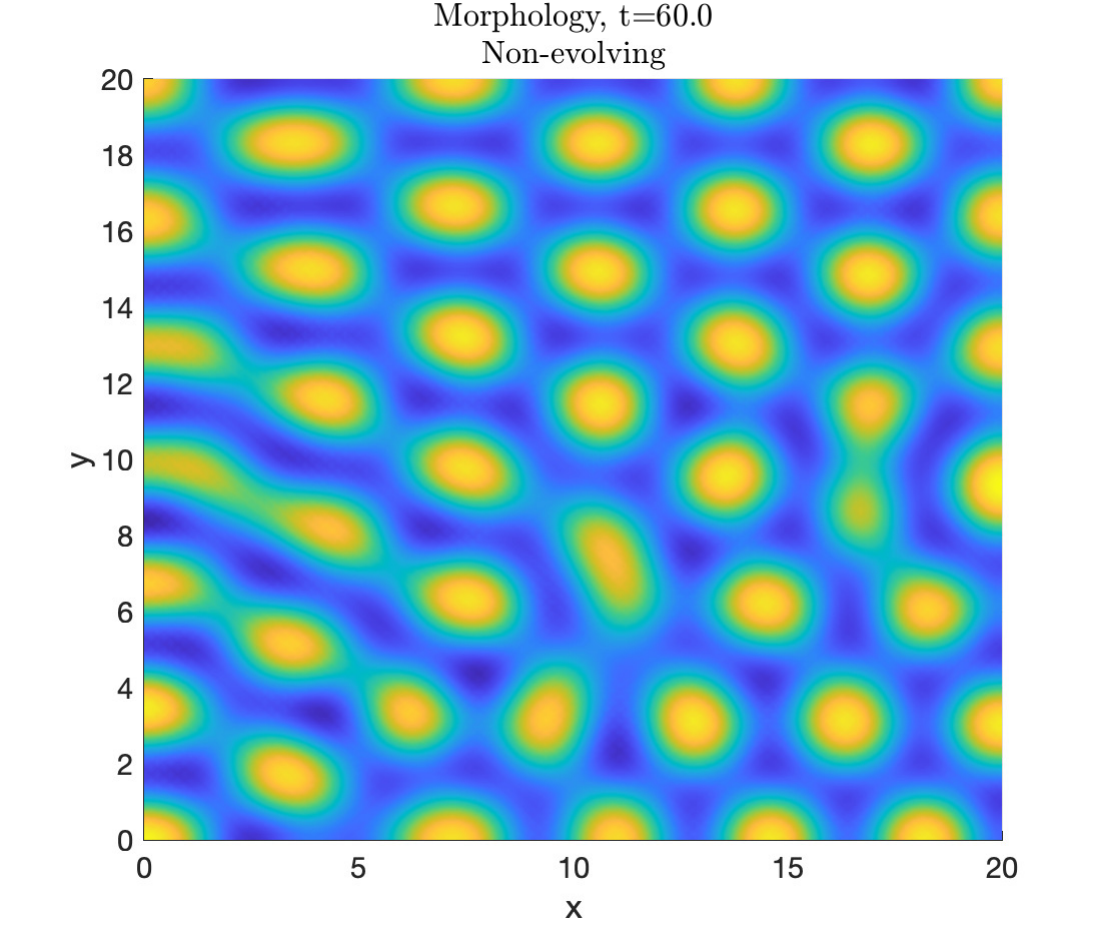}\hspace*{-5mm}
\includegraphics[scale=0.28]{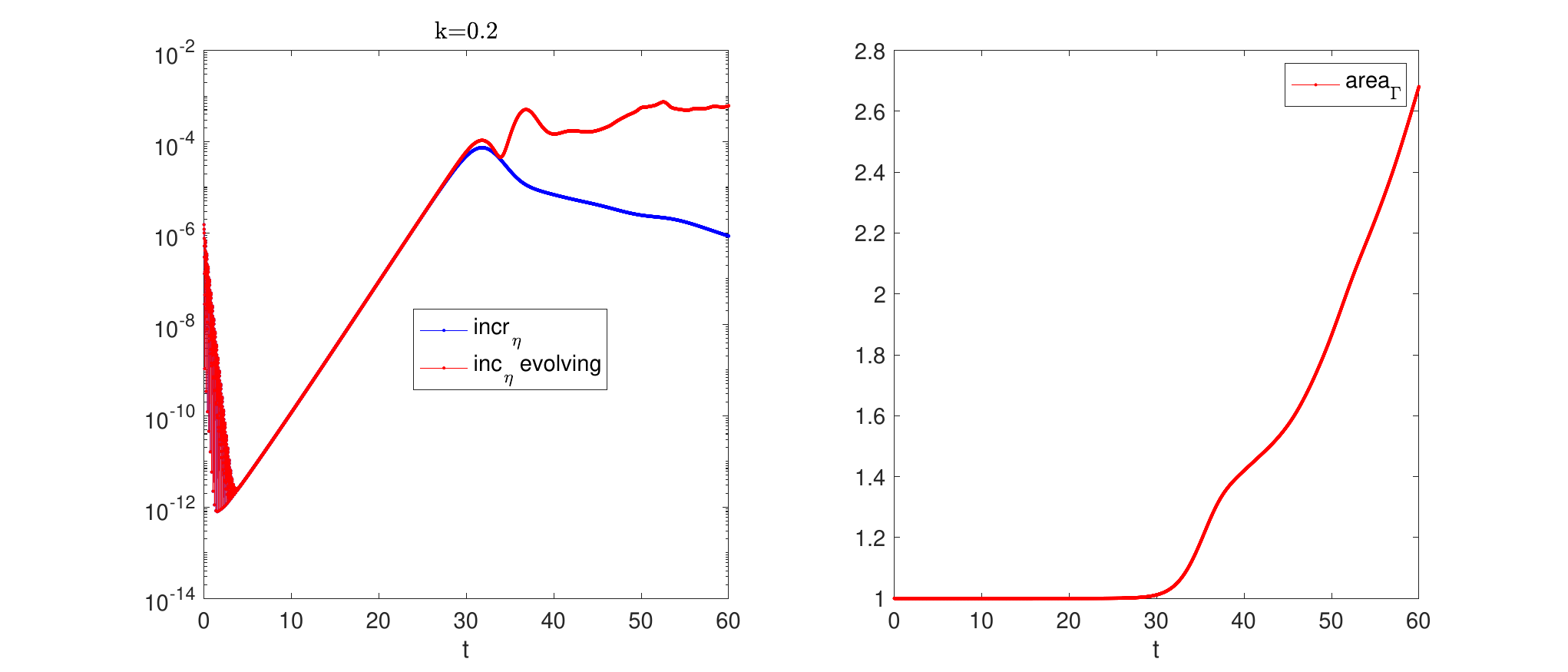}
\caption{Example 3. Left: $\eta$ component of the DIB model at the final time.  Middle: time increment of the numerical solutions of the DIB and ESDIB models.  Right: area of the evolving surface $\Gamma(t)$ in the ESDIB model.}
\label{fig:example_3_dib_and_increment}
\end{figure}

In this example, the initial domain is the square $\Omega = [0,L]^2$ for $L=20$. The final time is $T=60$. The \red{model} parameters are chosen as $B = 62$, $C=5$, \red{which lead to a holes-type pattern in the DIB model \cite{sgura2019parameter}}, and $\rho=2$. In Figure \ref{fig:example_3_esdib}, we show the numerical solution of the ESDIB model at different times (subfigure \ref{fig:example_3_esdib_num}) with the SEM micrograph of a zinc foil anode charged in a 6M KOH solution containing 100 ppm of tetra-butyl ammonium bromide at 5 mA cm$^{-2}$ for 1 hour (subfigure \ref{fig:example_3_esdib_lab}). This additive stabilizes the growth of compact layers of globular crystallites at intermediate-to-low current densities \cite{rossi2022insight}.  In Figure \ref{fig:example_3_dib_and_increment},  we show the solution of the DIB model at the final time, the time increments of the solutions of both models and the area of $\Gamma(t)$ in the ESDIB model over time.

\subsection{Example 4: sphere, $B=30$, $C=3$, holes}
\begin{figure}[p]
\begin{subfigure}{\textwidth}
\includegraphics[scale=0.35]{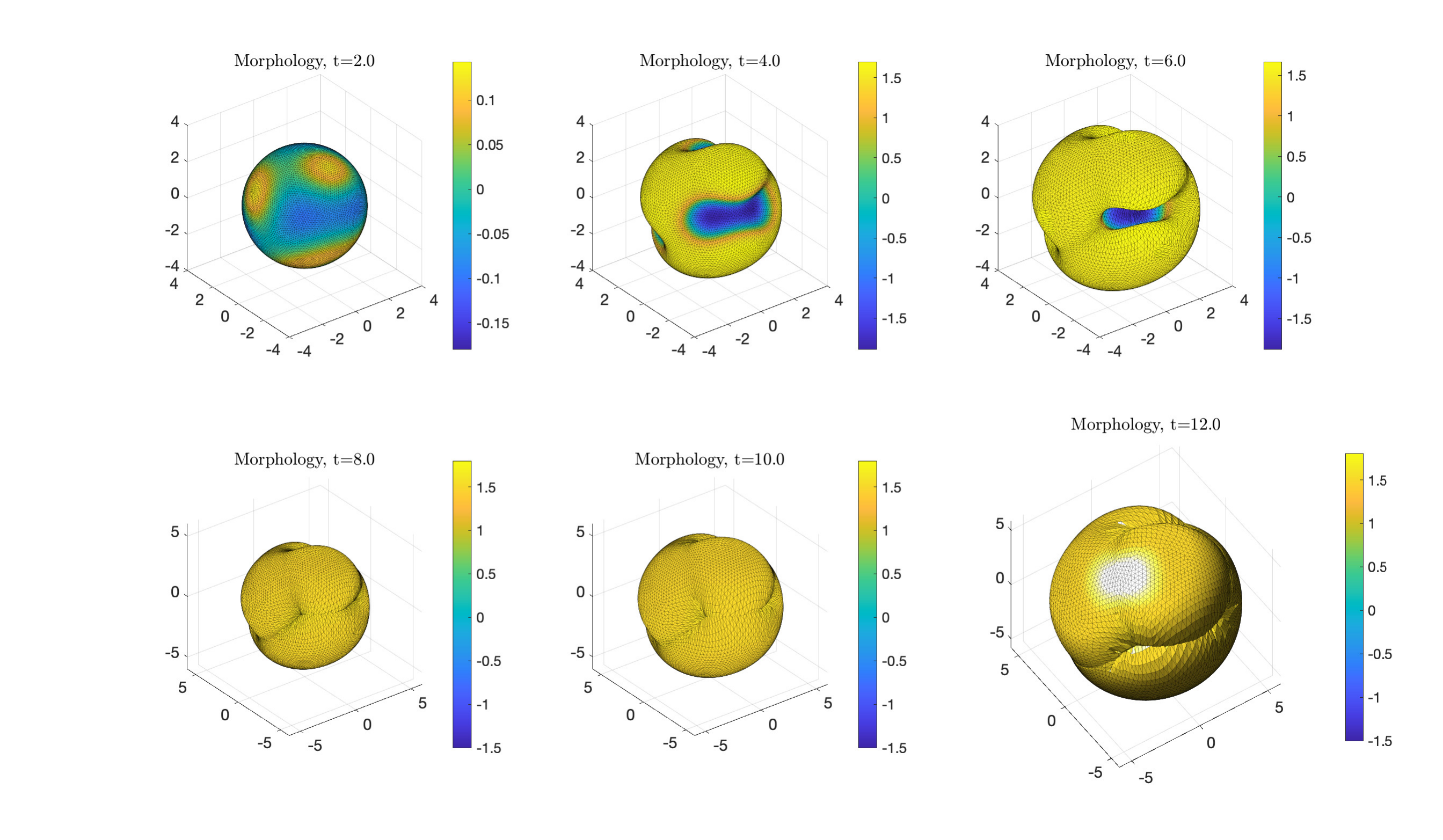}
\caption{}
\label{fig:example_4_esdib_num}
\end{subfigure}
\centering
\begin{subfigure}{\textwidth}
\centering
\includegraphics[scale=.25]{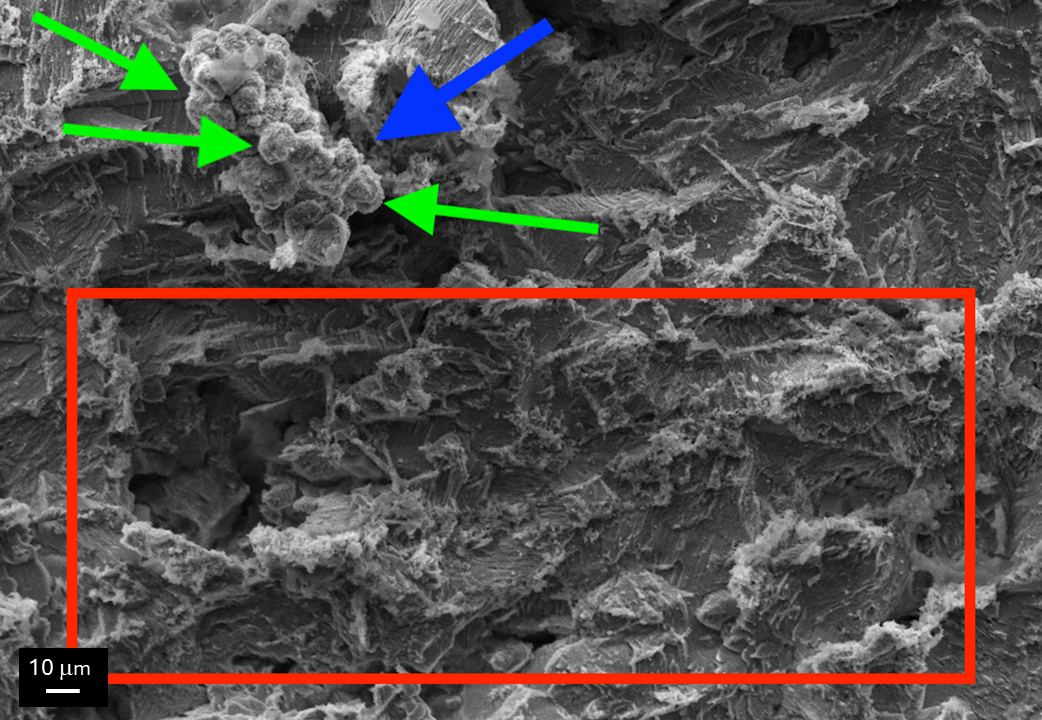}
\caption{}
\label{fig:example_4_esdib_lab}
\end{subfigure}
\caption{Example 4.  (a) $\eta$ component of the numerical solution of the ESDIB model at various times. (b) SEM micrograph of a zinc foil electrode subjected to ten discharge-charge cycles at 15 mA cm$^{-2}$ in pure 6M KOH and terminated with the charging cycle.}
\label{fig:example_4_esdib}
\end{figure}

\begin{figure}[p]
\centering
\includegraphics[scale=0.23]{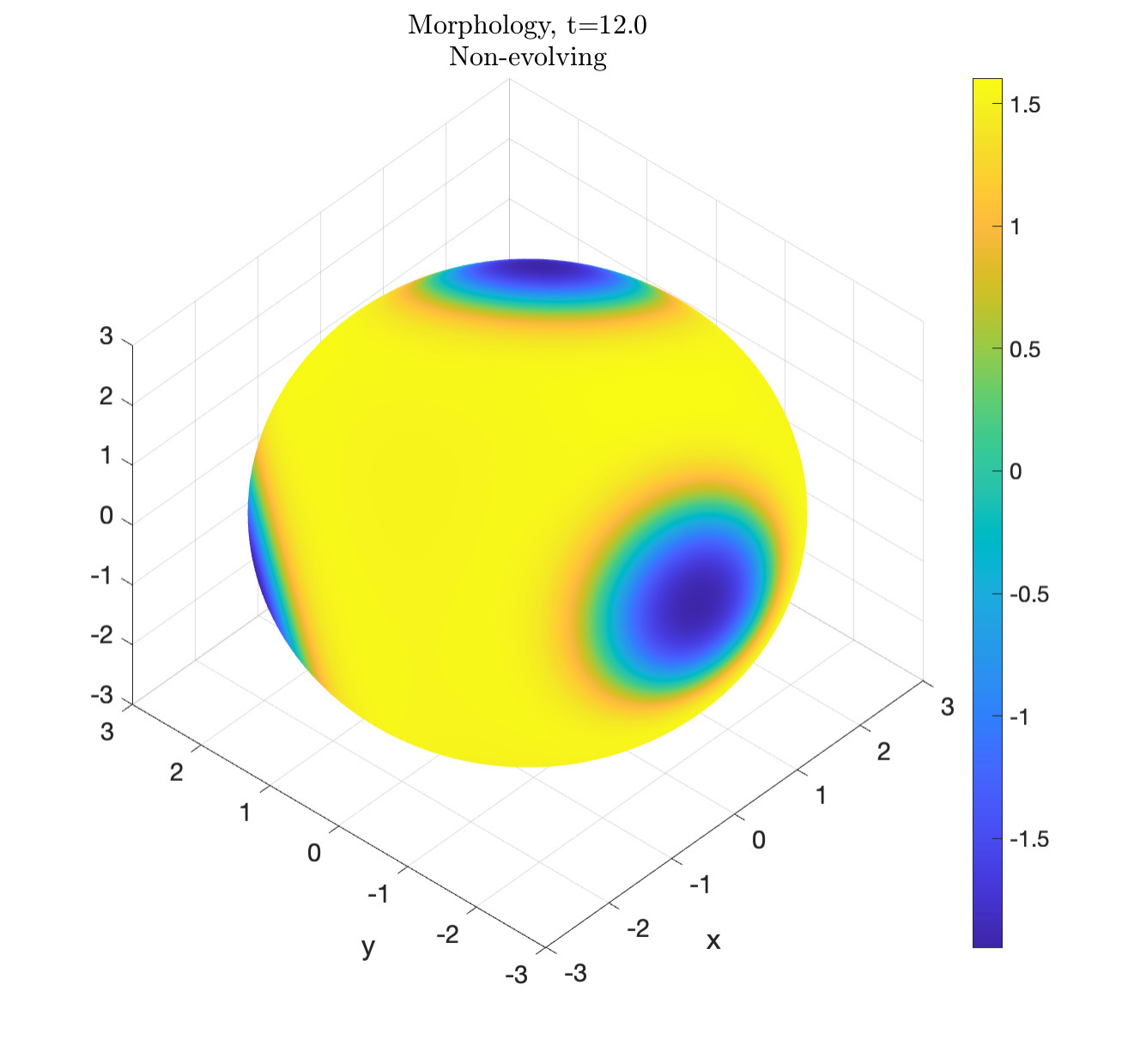}\hspace*{-5mm}
\includegraphics[scale=0.32]{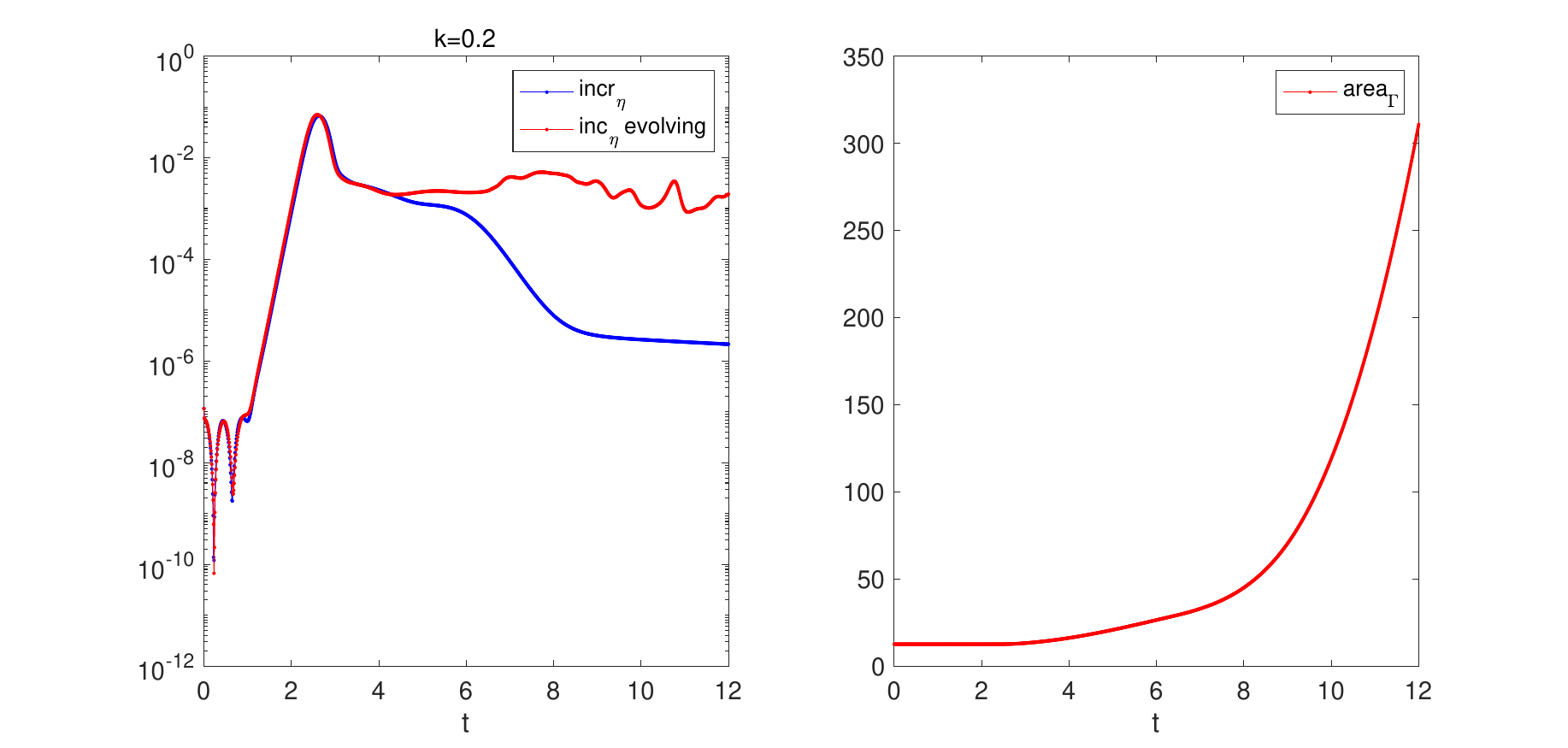}
\caption{Example 4.  Left: $\eta$ component of the DIB model at the final time.  Middle: time increment of the numerical solutions of the DIB and ESDIB models.  Right: area of the evolving surface $\Gamma(t)$ in the ESDIB model.}
\label{fig:example_4_dib_and_increment}
\end{figure}

In this example, the initial domain is the sphere $\Omega$ of radius for $R=3$. The final time is $T=12$. The \red{model}  parameters are chosen as $B = 30$, $C=3$, \red{which lead to a reversed spots (holes)-type pattern in the DIB model \cite{sgura2019parameter}}, and $\rho=2$.  In Figure \ref{fig:example_4_esdib}, we show the numerical solution of the ESDIB model at different times (subfigure \ref{fig:example_4_esdib_num}) and we report an SEM micrograph of a zinc foil electrode subjected to ten discharge-charge cycles at 15 mA cm$^{-2}$ in pure 6M KOH and terminated with the charging cycle (subfigure \ref{fig:example_4_esdib_lab}). These rather aggressive cycling conditions lead to localization of the plating and stripping processes \red{at} different positions of the electrode. Specifically (subfigure \ref{fig:example_4_esdib_lab}), extensive corroded regions can be noticed (red box), together with globular growth ones (blue arrow). Growth is hindered in the heavily corroded areas, owing to passivation, while renucleation occurs on the active globuli, leading to the formation of secondary globuli (thin, green arrows). For this reason, the growth process occurring during the last charging cycle is better captured by considering a spherical initial domain, as \red{illustrated} in subfigure \ref{fig:example_4_esdib_num}.  In Figure \ref{fig:example_4_dib_and_increment},  we show the solution of the SDIB model at the final time, the time increments of the solutions of both models and the area of $\Gamma(t)$ in the ESDIB model over time.

\subsection{Example 5: sphere, $B=66$, $C=3$, labyrinth}
\begin{figure}[p]
\begin{subfigure}{\textwidth}
\vspace*{-5mm}
\includegraphics[scale=0.35]{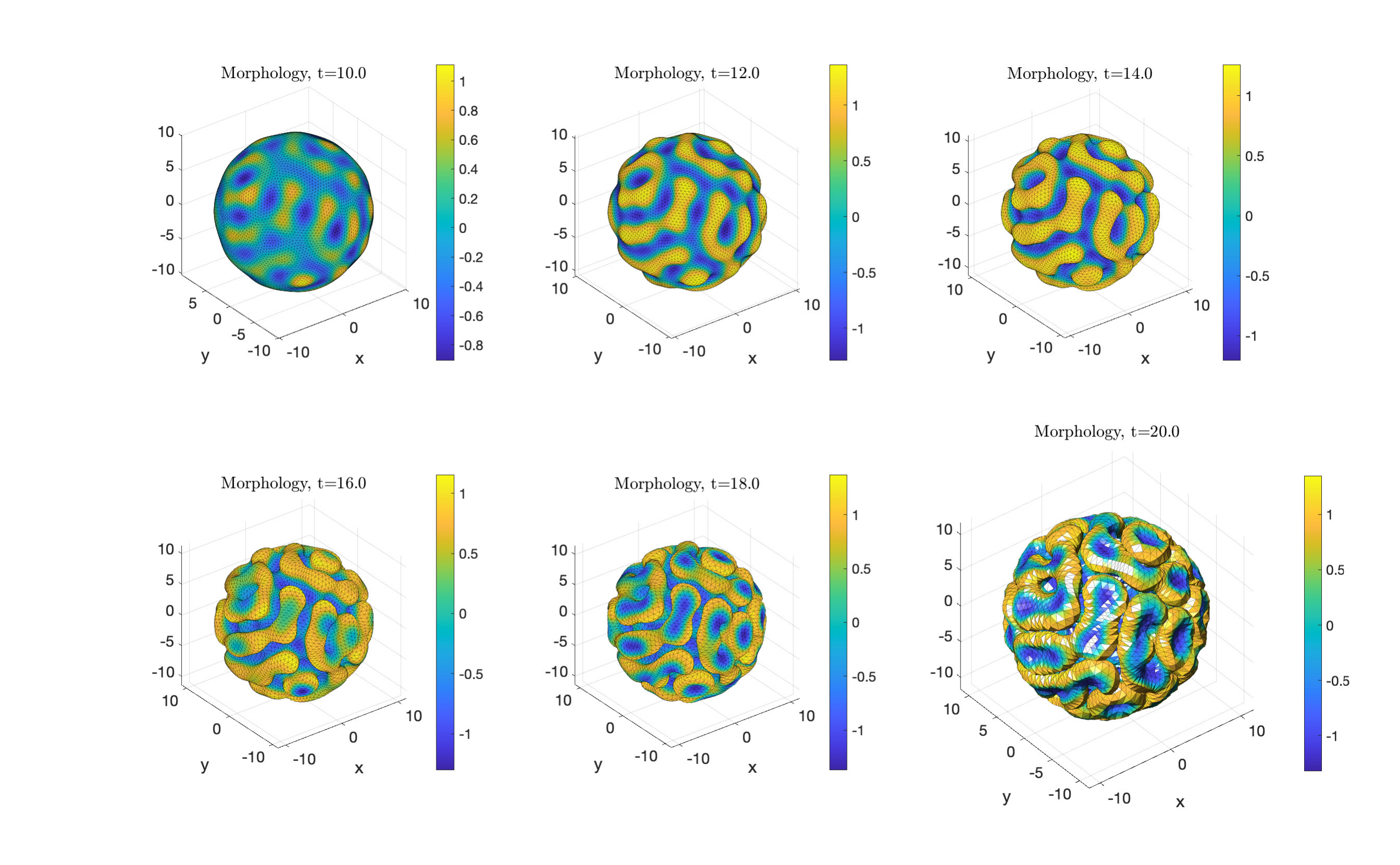}
\caption{}
\label{fig:example_5_esdib_num}
\end{subfigure}
\centering
\begin{subfigure}{\textwidth}
\centering
\includegraphics[scale=.25]{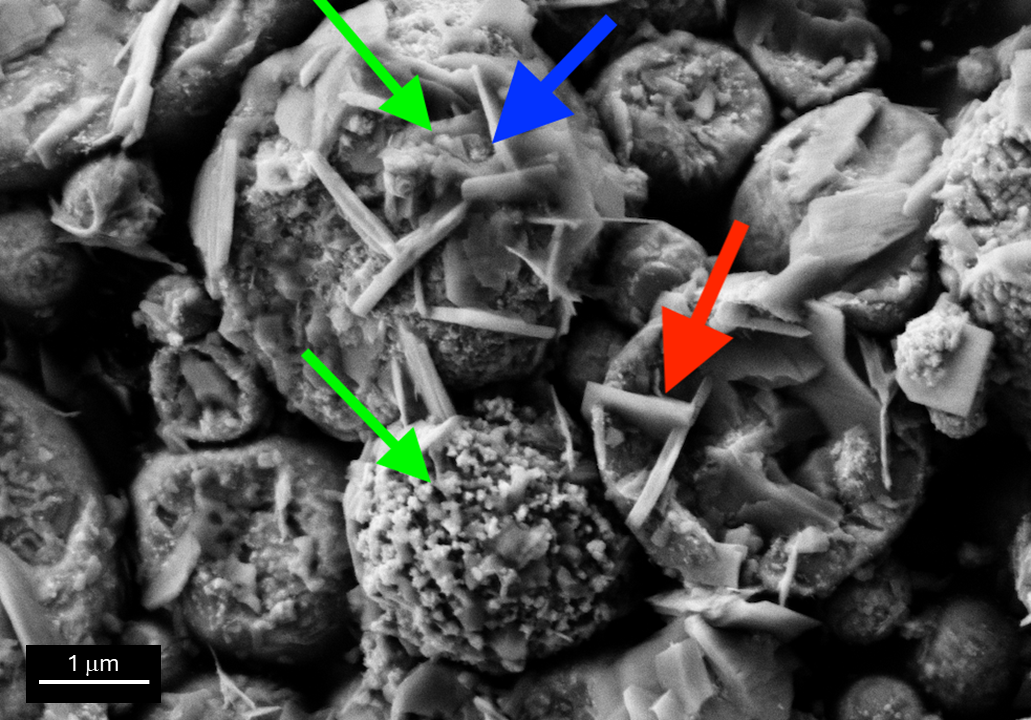}
\caption{}
\label{fig:example_5_esdib_lab}
\end{subfigure}
\caption{Example 5.  (a) $\eta$ component of the numerical solution of the ESDIB model at various times.  (b) SEM micrograph of a zinc sponge electrode subjected to 40 potentiostatic discharge-charge cycles of 1 hour at 50 mV vs. Zn in pure 6M KOH and terminated with charge.}
\label{fig:example_5_esdib}
\end{figure}

\begin{figure}[p]
\centering
\includegraphics[scale=0.27]{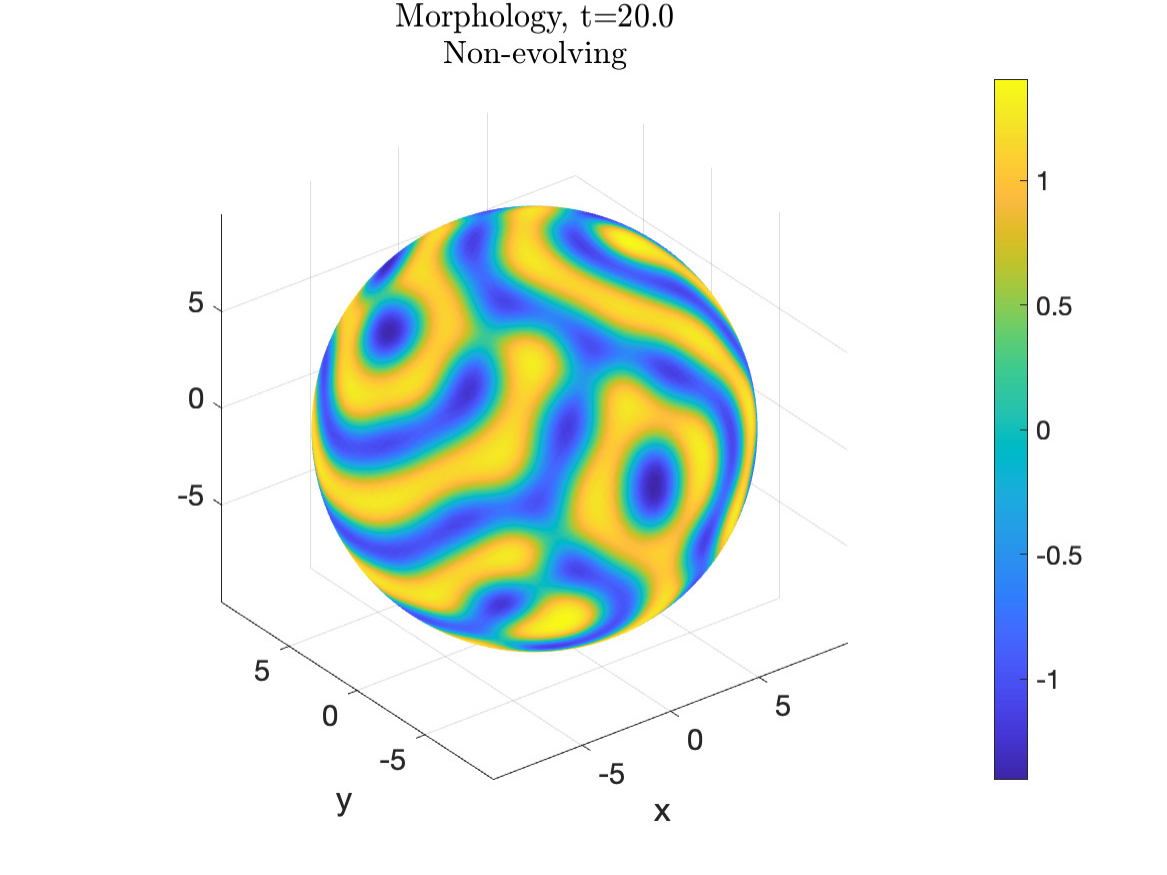}\hspace*{-5mm}
\includegraphics[scale=0.32]{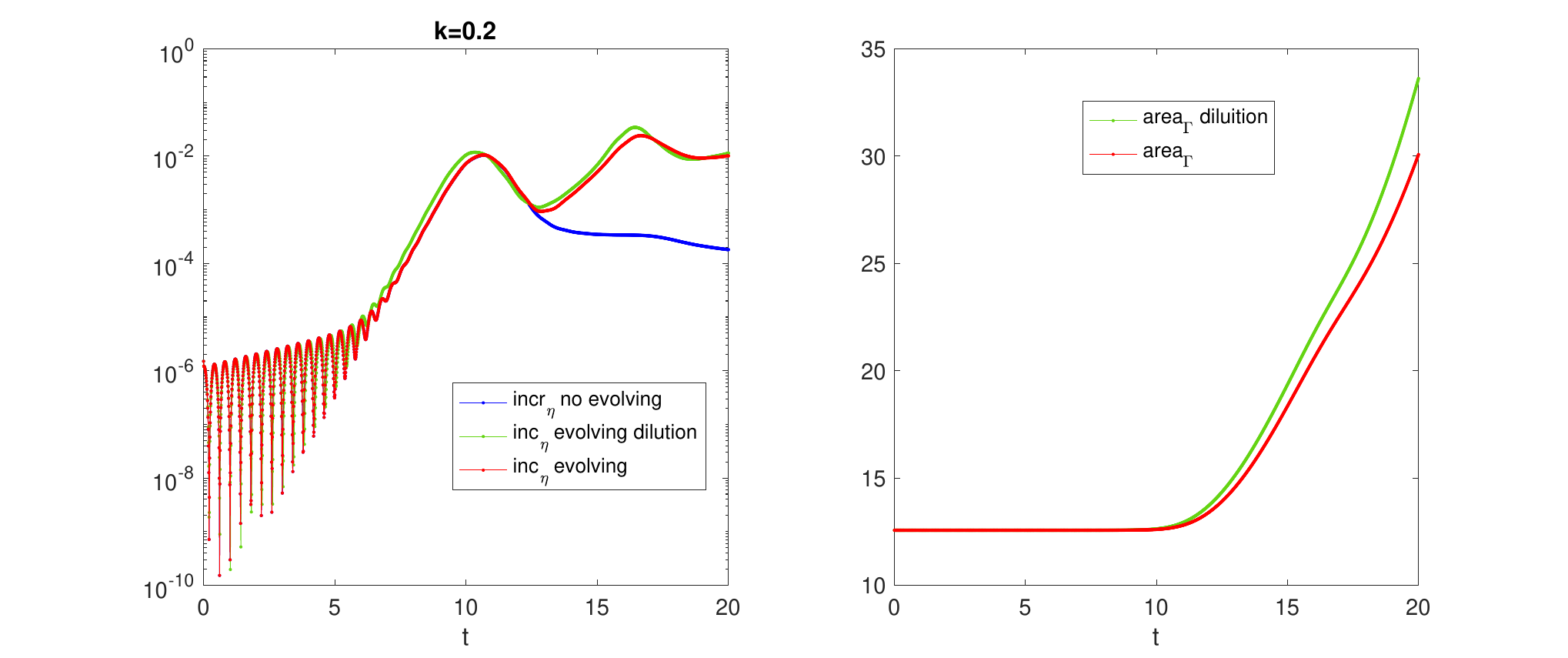}
\caption{Example 5.  Left: $\eta$ component of the DIB model at the final time.  Middle: time increment of the numerical solutions of the DIB and ESDIB models.  Right: area of the evolving surface $\Gamma(t)$ in the ESDIB model.}
\label{fig:example_5_dib_and_increment}
\end{figure}

In this example, the initial domain is the sphere $\Omega$ of radius for $R=10$. The final time is $T=20$. The \red{model}  parameters are chosen as $B = 66$, $C=3$, \red{which lead to a labyrinth-type pattern in the DIB model \cite{sgura2019parameter}}, and $\rho=1$. In Figure \ref{fig:example_5_esdib_num}, we show the numerical solution of the ESDIB model at different times. Figure \ref{fig:example_5_esdib_lab} reports an SEM micrograph of a zinc sponge anode subjected to 40 potentiostatic discharge-charge cycles of 1 hour at 50 mV vs.  Zn in pure 6M KOH and terminated with the charging cycle.  Specifically, charge-discharge cycling leads to templating with inactive regions --formed during the anodic interval--, characterized by ZnO crystals (red arrows) and active ones (blue arrows), where characteristic electrodeposition globuli form (thin, green arrows).  In Figure \ref{fig:example_5_dib_and_increment},  we show the solution of the SDIB model at the final time, the time increments of the solutions of both models and the area of $\Gamma(t)$ in the ESDIB model over time.

\subsection{Example 6: sphere, $B=62$, $C=5$, spots}
\begin{figure}[p]
\begin{subfigure}{\textwidth}
\includegraphics[scale=0.35]{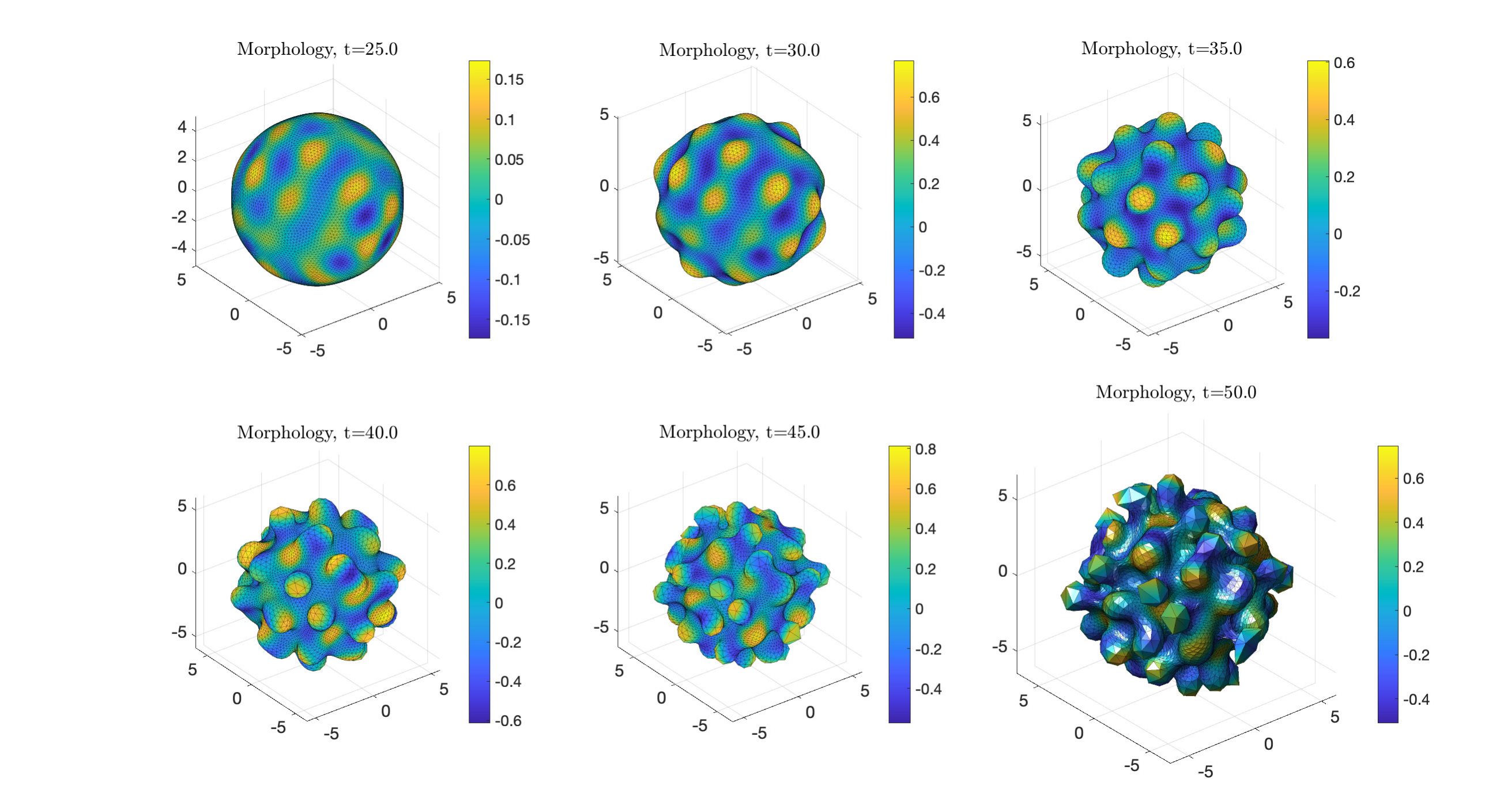}
\caption{}
\label{fig:example_6_esdib_num}
\end{subfigure}
\begin{subfigure}{\textwidth}
\centering
\includegraphics[scale=.25]{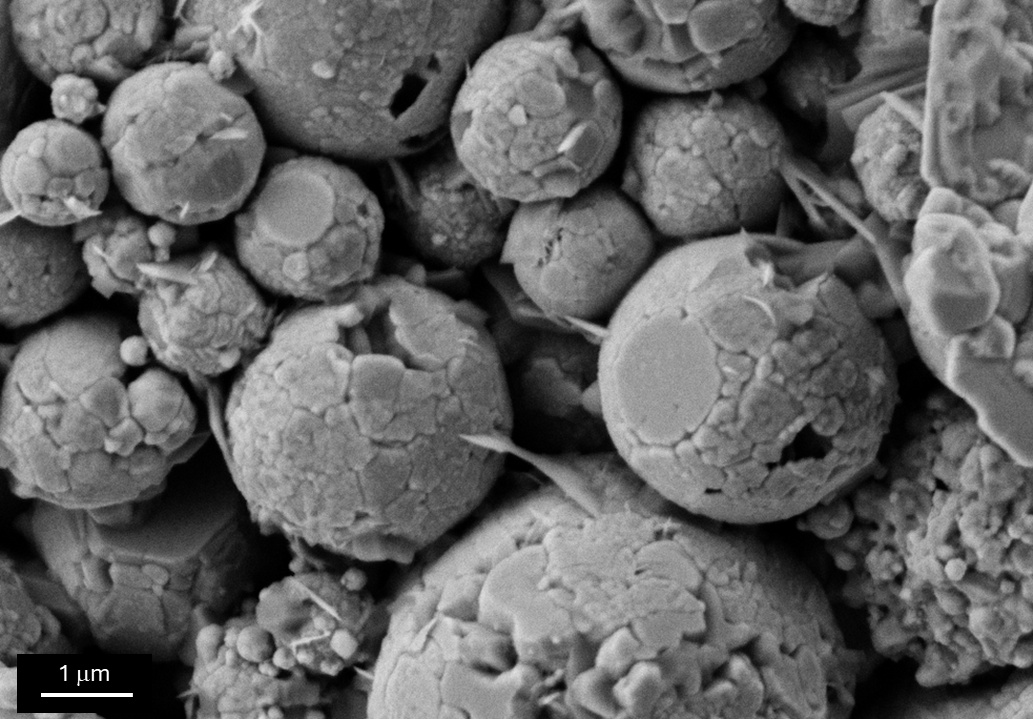}
\caption{}
\label{fig:example_6_esdib_lab}
\end{subfigure}
\caption{Example 6. (a) $\eta$ component of the numerical solution of the ESDIB model at various times.  (b) SEM micrograph of a zinc sponge electrode subjected first to a formation at -50 mV vs.  Zn for 2.5 h and then deep-discharged at 1300 mV vs.  Zn.}
\label{fig:example_6_esdib}
\end{figure}

\begin{figure}[p]
\centering
\includegraphics[scale=0.27]{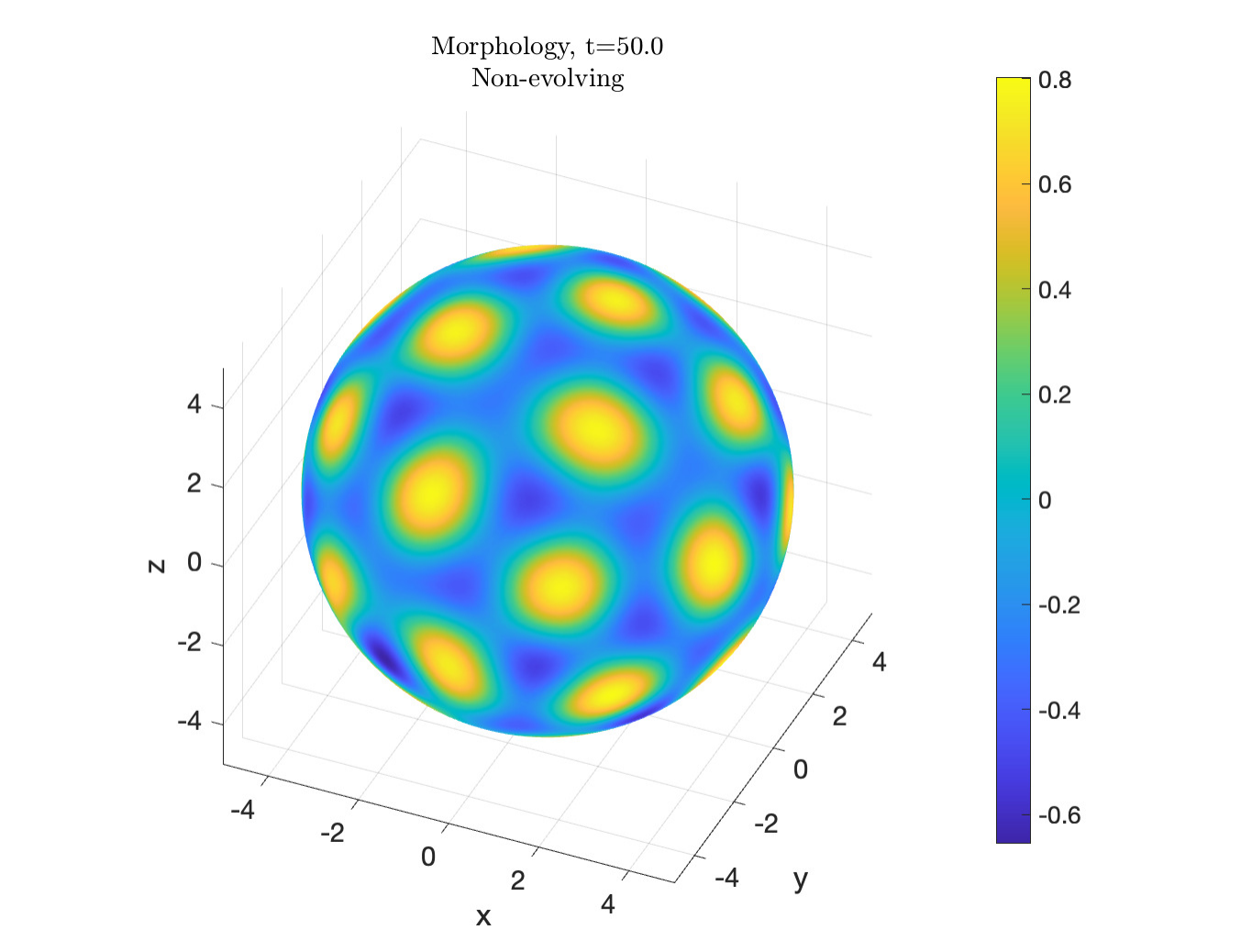}\hspace*{-5mm}
\includegraphics[scale=0.32]{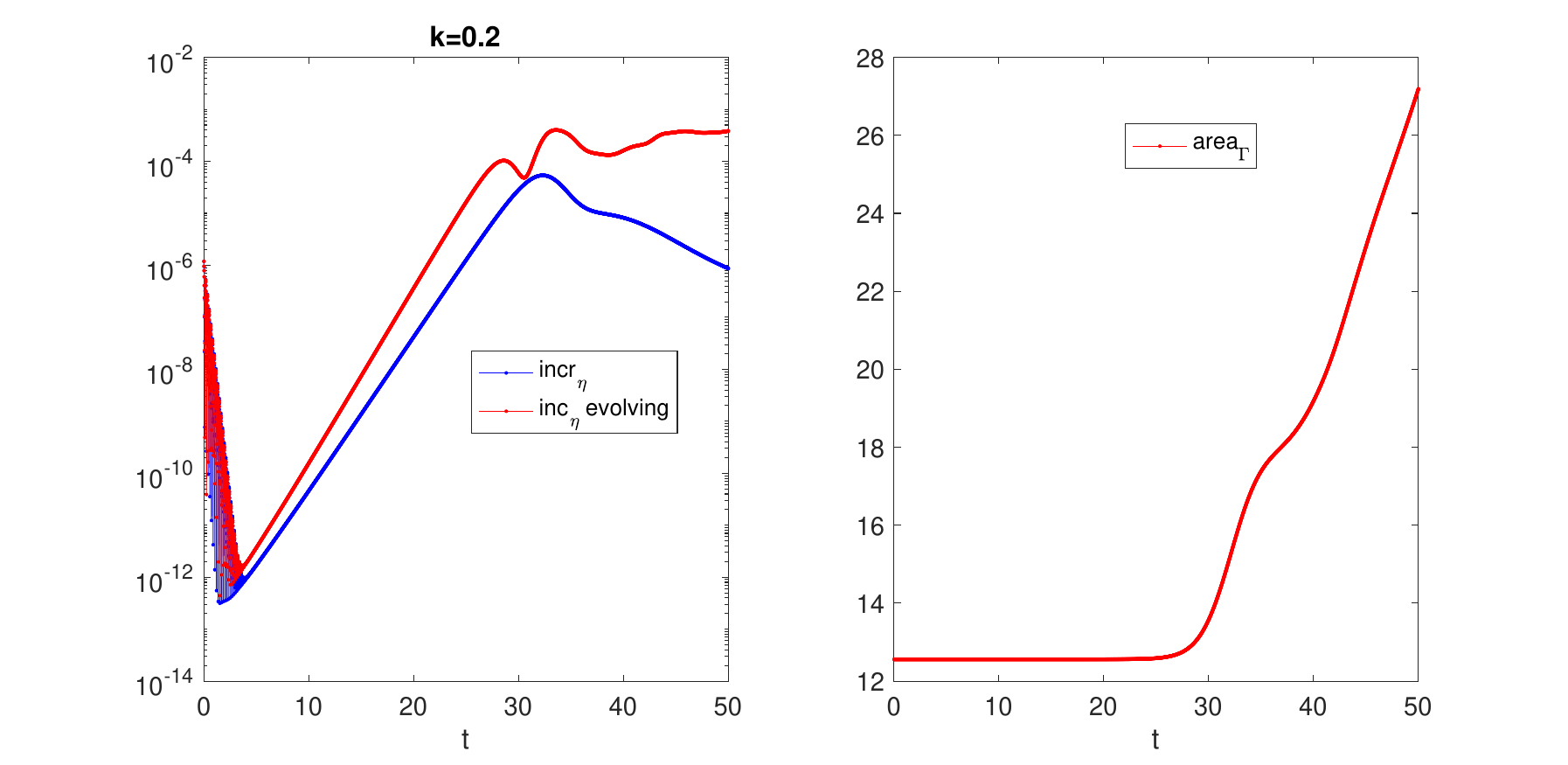}
\caption{Example 6.  Left: $\eta$ component of the DIB model at the final time.  Middle: time increment of the numerical solutions of the DIB and ESDIB models.  Right: area of the evolving surface $\Gamma(t)$ in the ESDIB model.}
\label{fig:example_6_dib_and_increment}
\end{figure}

In this example, the initial domain is the sphere $\Omega$ of radius for $R=5$. The final time is $T=50$. The \red{model}  parameters are chosen as $B = 62$, $C=5$, \red{which lead to a spots-type pattern in the DIB model \cite{sgura2019parameter}}, and $\rho=2$.  In Figure \ref{fig:example_6_esdib_num}, we show the numerical solution of the ESDIB model at different times.  Figure \ref{fig:example_6_esdib_lab} shows an SEM micrograph of a zinc sponge anode subjected first to formation at -50 mV vs. Zn for 2.5 h and then deep-discharged at 1300 mV vs.  Zn till zero-current conditions are virtually attained. The formation cycle leads to essentially full conversion of the pristine ZnO layer to metallic zinc with minor shape changes, and subsequent deep discharge brings about the conversion of metallic zinc to a thick ZnO layer. Nucleation and lattice mismatch jointly cause the formation of a granular crystallite structure onto the spheroidal electrode. In Figure \ref{fig:example_6_dib_and_increment},  we show the solution of the SDIB model at the final time, the time increments of the solutions of both models and the area of $\Gamma(t)$ in the ESDIB model over time.  

\FloatBarrier

\section{Conclusions}
\label{sec:conclusions}
In this work, we have introduced the \textit{Evolving Surface DIB (ESDIB)} model, a morphochemical reaction--diffusion system posed on a dynamically evolving electrode surface for battery modeling, to simulate metal electrodeposition and the formation of branching structures. By coupling unknown, morphology-driven surface evolution to morphology itself and local concentration of electrochemical species, the ESDIB framework overcomes the limitations of fixed-surface formulations, enabling a physically consistent description of morphological instabilities such as dendritic growth.  

Numerical simulations conducted with the \textit{Lumped Evolving Surface Finite Element Method (LESFEM)} for spatial discretisation and the \textit{IMEX Euler scheme} for time integration demonstrate good agreement with experimental microscopic images. The model successfully captures key qualitative and quantitative features of metal deposition, including the onset of branching and shape evolution, confirming its predictive capability and relevance for studying electrodeposition phenomena in energy storage devices.  We conclude that the interplay between mathematical modeling, numerical methods and  experimental evidence is key to advance the knowledge on material localization processes in batteries toward the design of more durable batteries.

A potential limitation of the current ESDIB model is its inability to correctly handle surface self-intersections, which can occur during extensive branching. In future work, we plan to address this limitation in two complementary ways. First, by developing an algorithm capable of managing topological changes in the evolving surface.  Second, by exploring an alternative \textit{phase-field type formulation}, which naturally accommodates complex interface dynamics and topological transitions.

\bibliographystyle{plainurl} 
\bibliography{bibliography.bib}

\end{document}